\documentclass[mathpazo]{cicp}

\usepackage{graphicx}
\usepackage{subfigure}
\usepackage{amssymb}
\usepackage{amsmath}
\usepackage{amsfonts}
\usepackage{mathrsfs}
\usepackage{mathtools}
\usepackage{empheq}
\usepackage{bm}
\usepackage{bbm}
\usepackage{color}
\usepackage{setspace}
\usepackage{epstopdf}
\usepackage{float}
\usepackage{cite}
\usepackage{nicefrac}
\usepackage{booktabs} 
\usepackage{array} 
\usepackage{paralist} 
\usepackage{verbatim} 
\usepackage{subfigure} 
\usepackage[normalem]{ulem}
\usepackage{stmaryrd}
\usepackage{multirow}
\usepackage{makecell}
\usepackage[colorlinks=true]{hyperref}

\usepackage{algorithm}
\usepackage{algorithmic}  

\usepackage{amsopn}
\DeclareMathOperator{\diag}{diag}

\def\softd{{\leavevmode\setbox1=\hbox{d}%
          \hbox to 1.05\wd1{d\kern-0.4ex{\char039}\hss}}}
          
\newcommand{\dpar}[2]{\dfrac{\partial #1}{\partial #2}}
\newcommand{\R}{\mathbb R}
\newcommand{\Z}{\mathbb Z}
\newcommand\eref[1]{(\ref{#1})}

\newcommand{\jph}{{j+\frac{1}{2}}}
\newcommand{\jmh}{{j-\frac{1}{2}}}
\newcommand{\dx}{\Delta x}

\newcommand{\dt}{\Delta t}

\newcommand*\xbar[1]{%
  \hbox{%
    \vbox{%
      \hrule height 0.5pt 
      \kern0.4ex
      \hbox{%
        \kern-0.05em
        \ensuremath{#1}%
        \kern-0.00em
      }%
    }%
  }%
}

\begin{document}
\title{Bound-preserving Point-Average-Moment PolynomiAl-interpreted (PAMPA) scheme: One-dimensional case}


 \author[R. Abgrall et.~al.]{R\'{e}mi Abgrall\affil{1},
       Miaosen Jiao\affil{2}, Yongle Liu \affil{1}\comma\corrauth, and Kailiang Wu\affil{3}}
 \address{\affilnum{1}\ Institute of Mathematics,
            University of Z\"{u}rich, 8057 Z\"{u}rich, Switzerland. \\
           \affilnum{2}\ Department of Mathematics, 
           Southern University of Science and Technology, Shenzhen 518055, China.\\
           \affilnum{3}\ Department of Mathematics and Shenzhen International Center for Mathematics, 
           Southern University of Science and Technology, Shenzhen 518055, China.}
 \emails{{\tt remi.abgrall@math.uzh.ch} (R.~Abgrall), {\tt 12332859@mail.sustech.edu.cn} (M.~Jiao),
          {\tt yongle.liu@math.uzh.ch; liuyl2017@mail.sustech.edu.cn} (Y.~Liu), {\tt wukl@sustech.edu.cn} (K.~Wu).}


\begin{abstract}
We propose a bound-preserving (BP) Point-Average-Moment PolynomiAl-interpreted (PAMPA) scheme by blending third-order and first-order constructions. 
The originality of the present construction is that it does not need any explicit reconstruction within each element, and therefore the construction is very flexible. The scheme employs a classical blending approach between a first-order BP scheme and a high-order scheme that does not inherently preserve bounds. 
The proposed BP PAMPA scheme demonstrates effectiveness across a range of problems, from scalar cases to systems such as the Euler equations of gas dynamics. We derive optimal blending parameters for both scalar and system cases, with the latter based on the recent geometric quasi-linearization (GQL) framework of [Wu \& Shu, {\em SIAM Review}, 65 (2023), pp. 1031--1073].  
This yields explicit, optimal blending coefficients that ensure positivity and control spurious oscillations in both point values and cell averages. 
This framework incorporates a convex blending of fluxes and residuals from both high-order and first-order updates, facilitating a rigorous BP property analysis. Sufficient conditions for the BP property are established, ensuring robustness while preserving high-order accuracy. Numerical tests confirm the effectiveness of the BP PAMPA scheme on several challenging problems.
\end{abstract}


\maketitle

\section{Introduction}

We consider the one-dimensional hyperbolic partial differential equations (PDEs) of conservation laws: 
\begin{equation}
\label{1.1}
\dpar{\mathbf u}{t}+\dpar{\mathbf f(\mathbf u)}{x}={\mathbf 0}
\end{equation}
with the initial condition $\mathbf u(x,0)=\mathbf u_0(x)$. Here, $x \in \Omega \subset \mathbb R$ is the spatial variable, $t \geq 0$ represents time, $\mathbf u(x,t) \in \mathbb{R}^n$ is an unknown vector function, and $\mathbf f(\mathbf u)$ denotes the physical fluxes.

It is known that the {entropy weak} solution \footnote{{An entropy weak solution should satisfy the following inequality in the distributional sense $\eta(\mathbf u)_t+q(\mathbf u)_x\leq0$, where $(\eta, q)$ is the so-called entropy pair.}} of \eref{1.1} often satisfies certain physical bounds. 
For instance, when \eqref{1.1} reduces to a scalar form, $u_t + f(u)_x = 0$, the entropy solution of scalar conservation laws adheres to the maximum principle. Specifically, if $u_0(x) \in [{\mathring{u}_{\min}, \mathring u_{\max}}]$, the entropy solution $u(x,t)$ will maintain these bounds for any $t > 0$, remaining within the invariant domain $\mathcal{D} = [{\mathring{u}_{\min}, \mathring u_{\max}}]$. In the case of a system of conservation laws, where $\mathbf u(x, t) \in \mathbb{R}^p$, one can often identify a convex invariant domain $\mathcal{D} \subset \mathbb{R}^p$ such that both $\mathbf u_0$ and $\mathbf u(x, t)$ lie within $\mathcal{D}$. 
In this work, we will focus on the one-dimensional compressible Euler equations, with
\begin{equation*}
  \mathbf u = \begin{pmatrix} \rho\\ m\\ E \end{pmatrix},
\end{equation*}
where $\rho$ represents the fluid density, $m = \rho v$ is the momentum, $v$ denotes the fluid velocity, and $E = e + \tfrac{1}{2} \rho v^2$ is the total energy, with $e$ as the internal energy. The flux vector is given by
\begin{equation*}
  \mathbf f(\mathbf u) = \begin{pmatrix} m \\ \frac{m^2}{\rho} + p \\ \frac{m (E + p)}{\rho} \end{pmatrix},
\end{equation*}
where $p = p(e, \rho)$ is the pressure, related to internal energy and density by an equation of state. For simplicity, we consider the ideal gas equation of state,
\begin{equation*}
    p = (\gamma - 1)e,
\end{equation*}
where $\gamma$ is a constant representing the specific heat ratio. For the Euler equations, a convex invariant domain $\mathcal{D}$ is defined by
\begin{equation}\label{Euler_IDP}
  \mathcal{D} = \left\{(\rho, m, E) \, \Big| \, \rho > 0 \text{ and } e = E - \frac{m^2}{2\rho} > 0 \right\}.
\end{equation}

It is highly desirable, and often crucial, to ensure that the numerical solution remains within the convex invariant domain $\mathcal{D}$. Numerical methods capable of preserving this convex set are referred to as bound-preserving (BP) or invariant-domain-preserving (or positivity-preserving in some special cases). Extensive work on this subject has been done in recent years; see, for instance, \cite{Lohmann_FCT,Zhang_MP,Carlier_IDP,Guermond_IDP3,Guermond_IDP,Guermond_IDP2,Hajduk_MCL,Kuzmin_MCL1,Kuzmin_BP,wu2023geometric,Vilar_DGFV}. Most of these cited works share a similar philosophy, namely blending high-order and low-order fluxes, operators, or schemes to ensure the preservation of convex properties and, more generally, to achieve BP properties. This is possible because certain first-order schemes have been shown to be BP for many hyperbolic systems. A suitable convex combination of first-order BP fluxes, operators, or schemes with high-order ones can be designed to achieve both the BP property and high-order accuracy simultaneously. 

Recently, a new class of schemes, capable of combining several formulations of the same hyperbolic problem \eqref{1.1} and inspired by the so-called Active Flux method (see, for instance, \cite{AF1,AF2,AF3,AF4,HKS,CHK_21}), was introduced in \cite{Abgrall_camc}. In this approach, the solution is represented continuously through the degrees of freedom (DoFs), which consist of average values within each element and point values at the element boundaries. These DoFs are then evolved according to two formulations of the same PDEs. The time stepping is simplified since the method of lines becomes possible. This novel approach has been successfully extended to the multidimensional case in \cite{Abgrall_Lin_Liu} and further applied to one-dimensional hyperbolic balance laws in \cite{AbgrallLiu, LiuBarsukow} and multidimensional hyperbolic balance laws in \cite{Liu_PAMPA_SW2D}, where this new procedure is named the 
PAMPA (Point-Average-Moment PolynomiAl-interpreted) scheme. 
Similar to high-order discontinuous Galerkin schemes, the third-order PAMPA scheme cannot prevent the numerical solution from producing spurious oscillations in the presence of discontinuities, which may lead to non-physical solutions, such as negative density and internal energy in gas dynamics. This issue can result in nonlinear instability and even code crashes. To address these problems, a stabilization technique equipped with parachute low-order schemes and a MOOD mechanism \cite{CDL} was applied in the aforementioned works. In \cite{duan2024AF}, a related problem was studied. The authors of \cite{duan2024AF} referred to the PAMPA scheme as a generalized Active Flux scheme and introduced two main contributions. First, the authors of \cite{duan2024AF} used the local Lax--Friedrichs, Steger--Warming, and van Leer--H\"{a}nel flux vector splittings to derive a high-order scheme for updating point values, addressing transonic issues in nonlinear problems caused by inaccurate upwind direction estimation in Jacobian splitting. Second, a convex limiting method is applied to the average DoFs, and a scaling limiter is used for the point value DoFs to ensure the BP property. In \cite{duan2024AF}, the low-order scheme for the point value is constructed using a first-order method based on the standard local Lax–Friedrichs scheme on a staggered mesh.

Observing from previous works on the PAMPA scheme (\cite{Abgrall_camc, Abgrall_Lin_Liu, AbgrallLiu, LiuBarsukow,Liu_PAMPA_SW2D}), we know that the parachute first-order schemes for both average and point value DoFs are BP under certain CFL condition. The aim of this paper is, therefore, to introduce a generic convex limiting framework for analyzing and constructing a BP PAMPA scheme, in which the high-order and low-order fluxes and residuals corresponding to average and point value DoFs will be blended using the same techniques. The efforts and findings include the following:
\begin{itemize}
    \item We propose a convex blending of the fluxes and residuals computed by the high-order PAMPA and low-order local Lax--Friedrichs' flavor schemes from \cite{Abgrall_camc}. With this blending, the updated solutions, including both average and point values, can be rewritten as a convex combination of the previous time step solution and several intermediate solution states. This formulation facilitates the BP property analysis. 
 {It is worth noting that, when updating the point value at the cell interface, our low-order scheme employs a robust local Lax–Friedrichs type scheme \cite{Abgrall_camc} that incorporates both point and average values, rather than the local Lax–Friedrichs scheme in \cite{duan2024AF}, which relies solely on point values. Additionally, our BP technique is based on a novel blending of residuals. Importantly, our BP PAMPA scheme does not exhibit the transonic issues seen in nonlinear problems in \cite{duan2024AF,HKS}, which required careful flux vector splitting in their approach in \cite{duan2024AF}.} 
    
    \item We establish a sufficient condition for the BP property of the blended PAMPA scheme, which only requires the BP property of the intermediate solution states. This finding converts the goal of constructing a BP PAMPA scheme into determining, through analysis, the optimal blending coefficients to reach the desired properties while maintaining as much of the high-order accuracy of the PAMPA scheme as possible. 
    \item We provide a set of theoretical results that can be used to rigorously prove the BP property of the BP PAMPA scheme. {Among these results, the analysis for maintaining the average value within the invariant domain aligns with that in \cite{duan2024AF}, while our BP limiting for the point value updates is different from \cite{duan2024AF}. Specifically, we write the scheme for the point value DoFs as a convex combination of two parts (residuals). We enforce that for each part the BP property is satisfied. Whereas in \cite{duan2024AF}, the authors used BP limiting on the whole solution state without a convex splitting.}

    \item We extend the BP PAMPA scheme to the Euler equations of gas dynamics. { Rather than using the techniques in \cite{duan2024AF, Vilar_DGFV, Kuzmin_MCL1},} we propose a novel approach based on the following property of $\mathcal{D}$ established in \cite{wu2023geometric}: setting 
 \begin{equation*}
     \mathcal{N}=\left\{\begin{pmatrix} 1\\ 0 \\ 0\end{pmatrix}\right\}\cup\left\{ \begin{pmatrix}\frac{\nu^2}{2} \\ -\nu\\ 1\end{pmatrix}, \nu\in \R\right\},
 \end{equation*}
the equivalent GQL (Geometric Quasi-Linearization) representation of the invariant domain $\mathcal{D}$ in \eref{Euler_IDP} is given by
\begin{equation}\label{Euler_GQL}
   \mathcal{D}_{\nu}=\{\mathbf u=(\rho, m, E)^T \text{ such that for all } n\in \mathcal{N}, \mathbf u^Tn>0\}. 
\end{equation}
The optimal blending coefficients can be directly obtained from the eigenvalues of a matrix with components given only by the intermediate solution states and the differences between high-order and low-order fluxes or residuals. This analysis is novel and applicable to the multidimensional case (which will be discussed in detail in a separate paper).
 \item  We implement our BP PAMPA scheme and demonstrate its robustness and effectiveness on several challenging numerical examples, such as the Buckley--Leverett problem in scalar conservation law with non-convex flux, the LeBlanc shock tube, double rarefaction Riemann, and Sedov problems in the Euler equations of gas dynamics.  
\end{itemize}

The structure of this paper is as follows. In Section \ref{sec2}, we review the one-dimensional first-order and high-order (here third-order) PAMPA schemes for solving \eqref{1.1}. We then use a blending procedure to combine the high-order and low-order schemes, demonstrating how to construct a BP PAMPA scheme. Specifically, we begin with scalar conservation laws, proposing BP blended fluxes for updating cell averages in Subsection \ref{subsec:average} and BP blended residuals for updating point values in Subsection \ref{subsec:point}. To reduce spurious oscillations near discontinuities while maintaining high accuracy in smooth regions, we introduce the local maximum principle and smooth extrema detector in Subsection \ref{subsec33}. The BP PAMPA scheme is then extended to the one-dimensional Euler equations of gas dynamics in Subsection \ref{EulerGQL}. Several challenging numerical examples are presented in Section \ref{sec4} to demonstrate the efficiency and robustness of the proposed BP PAMPA scheme. Finally, we conclude the paper in Section \ref{sec5} by summarizing the contributions and discussing future directions.

\section{One-dimensional PAMPA scheme for hyperbolic conservation laws}\label{sec2}
This section is devoted to recall the essential ingredients of the one-dimensional semi-discrete PAMPA scheme. To remain as simple as possible, the one-dimensional scalar conservation laws will be considered in this section. The system extension is perfectly straightforward. Spatial discretization of the studied PDE will be carried out using either a first-order scheme or a high-order scheme. To evolve the solution in time, we employ {the third-order strong stability-preserving (SSP) Runge--Kutta (RK) time-stepping method in the numerical implementation. However, for ease of demonstration, we apply only the first-order forward Euler method here.  Extending this to higher-order SSP RK time-stepping methods is conceptually straightforward, since these multistage time integration methods can be written as convex combinations of formal forward Euler steps.} 

For the subsequent discretization, let us introduce the following notations. The one-dimensional domain $\Omega$ is divided into a set of non-overlapping cells $I_\jph=[x_j, x_{j+1}]$, with $\dx_\jph=x_{j+1}-x_j$ being the size of $I_\jph$. We also partition the time domain in intermediate times $\{t^n\}$ with $\dt^n=t^{n+1}-t^n$ the $n$-th time step. Following \cite{Abgrall_camc,Abgrall_Lin_Liu}, the problem \eref{1.1} is discretized, at first order, by
\begin{subequations}
\label{scheme:O1}
\begin{equation}\label{scheme:O1a}
\xbar u_\jph^{n+1}=\xbar u_\jph^n-\lambda\Big ( \mathcal{F}^{LO}_{j+1}\big(\xbar u_\jph^n, \xbar{u}_{j+\frac{3}{2}}^n\big)-\mathcal{F}^{LO}_j\big(\xbar u_\jmh^n, \xbar u_\jph^n\big)\Big ),
\end{equation}
and
\begin{equation}\label{scheme:O1b}
u_j^{n+1}=u_j^n-2\lambda\Big(\overrightarrow{\Phi}_\jmh^{LO}+\overleftarrow{\Phi}_\jph^{LO}\Big).
\end{equation}
\end{subequations}
Here, $u_j^n$ is an approximation of the point value at $x_j$, $\xbar u_\jph^n$ is an approximation of the cell average in $[x_j,x_{j+1}]$, $\lambda=\frac{\Delta t^n}{\Delta x}$ (we assume a uniform mesh but this is not essential). In \eref{scheme:O1a}, the numerical flux $\mathcal{F}$ is monotone of Lipschitz constant $L$, and we make use of the well-known local Lax-Friedrichs numerical flux
\begin{equation}\label{LLF_flux}
\mathcal{F}^{LO}_j:=\mathcal{F}^{LO}_j\big(\xbar{u}_\jmh^n,\xbar{u}_\jph^n\big)=\frac{f(\xbar{u}_\jmh^n)+f(\xbar{u}_\jph^n)}{2}-\frac{\alpha_j}{2}\big(\xbar{u}_\jph^n-\xbar{u}_\jmh^n\big),
\end{equation}
where 
\begin{equation}\label{alpha_ave}
  \alpha_j=\alpha(\xbar{u}_\jmh^n,\xbar{u}_\jph^n)=\max\big(\vert f'(\xbar{u}_\jmh^n)\vert,\vert f'(\xbar{u}_\jph^n)\vert\big).
\end{equation}
In \eref{scheme:O1b}, the first-order residuals are computed by the following manner which also has the flavor of local Lax-Friedrichs scheme:
\begin{equation}\label{residu_LO}
  \left\{\begin{aligned}
  \overrightarrow{\Phi}_\jmh^{LO}&=\frac{f(u_j^n)-f(\xbar u_\jmh^n)}{2}-\frac{\beta_\jmh}{2}\big(\xbar u_\jmh^n-u_j^n\big),\\
   \overleftarrow{\Phi}_\jph^{LO}&=\frac{f(\xbar u_\jph^n)-f(u_j^n)}{2}-\frac{\beta_\jph}{2}\big(\xbar u_\jph^n-u_j^n\big),
  \end{aligned}\right.
\end{equation}
where 
\begin{equation}\label{beta_point}
  \beta_\jph=\beta(u_j^n,\xbar{u}_\jph^n)=\max\big(\vert f'(u_j^n)\vert,\vert f'(\xbar{u}_\jph^n)\vert\big).
\end{equation}
Under condition $\lambda\cdot \max\limits_j\alpha_j\leq 1/2$, $\lambda\cdot \max\limits_j\beta_\jph\leq 1/2$, and $L\lambda\leq 1$, the scheme is BP.

A high order version from \cite{Abgrall_camc} is
\begin{subequations}
\label{scheme:high}
\begin{equation}\label{scheme:high:average}
\xbar u_\jph^{n+1}=\xbar u_\jph^n-\lambda\Big ( f\big(u_{j+1}^n\big)-f\big(u_j^n\big)\Big )
\end{equation}
and
\begin{equation}\label{scheme:high:points}
u_j^{n+1}=u_j^n-2\lambda\Big(\overrightarrow{\Phi}_\jmh^{HO}+\overleftarrow{\Phi}_\jph^{HO}\Big)
\end{equation}
\end{subequations}
with the high-order residuals being defined as
\begin{equation}\label{HO_rd}
\overrightarrow{\Phi}_\jmh^{HO}=\frac{\dx}{2}a_j^+\dpar{u}{x}\big|_{ [x_{j-1},x_j]}(x_j),\quad
\overleftarrow{\Phi}_\jph^{HO}=\frac{\dx}{2}a_j^-\dpar{u}{x}\big|_{[x_{j},x_{j+1}]}(x_j),
\end{equation}
where $a_j=f'(u_j^n)$, and {$a_j^+=\max(a_j,0)$, $a_j^-=\min(a_j,0)$}. $u\big|_{[x_j,x_{j+1}]}$ is the approximation of $u$ in $[x_j,x_{j+1}]$:
\begin{equation}\label{u_para}
    u\big|_{[x_j,x_{j+1}]}=u_j^n \ell_0\Big ( \frac{x-x_j}{\Delta x}\Big ) +\xbar{u}_\jph^n \ell_{1/2}\Big ( \frac{x-x_j}{\Delta x}\Big )+u_{j+1}^n \ell_1\Big ( \frac{x-x_j}{\Delta x}\Big )
\end{equation}
with
\begin{equation*}
    \ell_0=(1-x)(1-3x), \quad \ell_{1/2}=6 x(1-x), \quad \ell_1(x)=x(3x-2).
\end{equation*}
These functions satisfy
\begin{equation*}
\begin{aligned}
&\ell_0(0)=1, && \ell_0(1)=0, && \int_0^1 \ell_0{\rm d}x=0, \\
&\ell_0(1)=0, && \ell_1(1)=1, && \int_0^1 \ell_1{\rm d}x=0, \\  
&\ell_{1/2}(0)=0,&& \ell_{1/2}(1)=0, &&\int_{0}^1\ell_{1/2}{\rm d}x=1.
\end{aligned}
\end{equation*}
Using the parabolic approximation in \eref{u_para}, we find
\begin{equation*}
\begin{aligned}
   &\dpar{u}{x}\big|_{[x_{j-1},x_j]}(x_j)=\frac{2}{\Delta  x}\big ( u_{j-1}^n-3\xbar u_\jmh^n+2u_j^n\big ), \\
    &\dpar{u}{x}\big|_{[x_{j},x_{j+1}]}(x_j)=\frac{2}{\Delta x}\big ( -2u_j^n+3\xbar{u}_\jph^n-u_{j+1}^n\big ).
\end{aligned} 
\end{equation*}
Then the high-order residuals in \eref{HO_rd} can be rewritten as
\begin{equation}\label{residuals_HO}
\left\{\begin{aligned}
\overrightarrow{\Phi}_\jmh^{HO}&=a_j^+(u_{j-1}^n-3\xbar u_\jmh^n+2u_j^n),\\
  \overleftarrow{\Phi}_\jph^{HO}&=a_j^-(-2u_{j}^n+3\xbar u_\jph^n-u_{j+1}^n).
\end{aligned}\right.
\end{equation}

In this procedure, the update of the conservative DoFs, $\xbar{u}_{j+1/2}$, must be done in a conservative manner. The conservation constraint can be relaxed for the update of the point values $u_j$. We have some freedom to select the point values at the cell boundaries, it can be either conservative variables, or primitive variables, or entropy variables as shown in \cite{Abgrall_camc}. It has also been demonstrated that under mild assumptions on the scheme, and assumptions similar to those of the Lax-Wendroff theorem on $\{u_j, \xbar u_{j+1/2}\}_{j\in \Z}$, the solution will converge to a weak solution of the problem \eref{1.1}. In this section as well as in this paper, we only consider the conservative variables for both the cell average and point value and the high-order scheme is the third-order one. However, the PAMPA method and the BP limitings developed in the subsequent section are not restricted by this. The extension to employ other variables for point value and higher-order case by using the higher moments will be investigated in the future work. 

When the procedure is extended to Euler equations of gas dynamics, we again apply the Jacobian splitting and replace $a_j^\pm$ by $\big(J(u_j)\big)^\pm=R\Lambda^\pm R^{-1}$, where $\Lambda=\diag(\lambda_1,\ldots,\lambda_p)$, $\Lambda^\pm=\diag(\lambda_1^\pm,\ldots,\lambda_p^\pm)$ {with $\lambda_i^+=\max(\lambda_i,0)$ and $\lambda_i^-=\min(\lambda_i,0)$}. $\lambda_1, \ldots, \lambda_p$ are the eigenvalues, with the columns of $R$ the corresponding eigenvectors, of the Jacobian matrix $J(u)=\frac{\partial f(u)}{\partial u}$.

In the numerical results, we will show that this procedure may lead to non-physical solutions. This has been observed by several authors \cite{HKS,duan2024AF}. However, we will show that our BP technique cures completely this problem.

\section{Bound-preserving PAMPA scheme}\label{sec3}
In this section, we introduce the BP PAMPA scheme. The previous presentation of the first-order and high-order schemes allows us to construct a unified blending framework for updating both cell average and point value. For simplicity and clarity in illustrating the blending procedure, the following subsections \ref{subsec:average}, \ref{subsec:point}, and \ref{subsec33} focus on scalar conservation laws. {For scalar conservation laws, if the initial data satisfies $ u_0\in\mathcal{D}=[\mathring{u}_{\min}, \mathring u_{\max}]$, then due to Kruzhkov's theory, the unique entropy weak solution $u(x,t)$ guarantees $u(\cdot,t)\in\mathcal{D}$, for any time $t$. This justifies why we want to construct schemes that preserves this property.} The extension to Euler equations will be addressed in subsection \ref{EulerGQL}.

\subsection{Blended fluxes and update of cell average}\label{subsec:average}
We present a convex limiting approach to achieve the BP property for updating the cell average, {inspired by \cite{Vilar_DGFV}}. To this end, each cell interface $x_j$ is assigned two fluxes: a first-order finite volume numerical flux $\mathcal{F}_j^{LO}(\xbar{u}_\jmh^n,\xbar{u}_\jph^n)$, and a high-order exact flux function $f(u_j^n)$. Then, these two fluxes are blended in a convex manner through a blending coefficient $\eta_j\in[0,1]$ as follows:
\begin{equation}\label{eq:blend_flux}
\mathcal{F}_j=\mathcal{F}_j^{LO}+\eta_j\big(f(u_j^n)-\mathcal{F}_j^{LO}\big)=\mathcal{F}_j^{LO}+\eta_j\Delta \mathcal{F}_j.
\end{equation}
A blending coefficient set to zero will lead to a first-order finite volume numerical flux, while a coefficient of one yields the high-order exact flux function. 

The first-order forward Euler time integration for the cell average reads as
\begin{equation}\label{eq:ave_new}
    \xbar{u}_\jph^{n+1}=\xbar{u}_\jph^n-\lambda(\mathcal{F}_{j+1}-\mathcal{F}_j),\quad \forall j,
\end{equation}
which can {be} rewritten as a convex combination of quantities defined at the previous time step: 
\begin{equation}\label{u_ave_new}
    \begin{aligned}
\xbar{u}_\jph^{n+1}&=\xbar{u}_\jph^n-\lambda(\mathcal{F}_{j+1}-\mathcal{F}_j)\pm\lambda(\alpha_j+\alpha_{j+1})\xbar{u}_\jph^n+\lambda\big(f(\xbar{u}_\jph^n)-f(\xbar{u}_\jph^n)\big)\\
&=(1-\lambda\alpha_j-\lambda\alpha_{j+1})\xbar{u}_\jph^n+\lambda\alpha_j\Big(\xbar{u}_\jph^n-\frac{-\mathcal{F}_j+f(\xbar{u}_\jph^n)}{\alpha_j}\Big)\\
&\quad +\lambda\alpha_{j+1}\Big(\xbar{u}_\jph^n-\frac{\mathcal{F}_{j+1}-f(\xbar{u}_\jph^n)}{\alpha_{j+1}}\Big)\\
&=(1-\lambda\alpha_j-\lambda\alpha_{j+1})\xbar{u}_\jph^n+\lambda\alpha_j{\widetilde{u}_{\jph}^L}+\lambda\alpha_{j+1}{\widetilde{u}_{\jph}^R},
    \end{aligned}
\end{equation}
with
\begin{equation*}
\begin{aligned}
  {\widetilde{u}_\jph^L}&=\xbar{u}_\jph^n-\frac{-\mathcal{F}_j+f(\xbar{u}_\jph^n)}{\alpha_j}
  =\frac{\xbar{u}_\jmh^n+\xbar{u}_\jph^n}{2}-\frac{f(\xbar{u}_\jph^n)-f(\xbar{u}_\jmh^n)}{2\alpha_j}+{\eta_j\frac{\Delta\mathcal{F}_j}{\alpha_j}}\\
  &=u_j^*+\eta_j\frac{\Delta\mathcal{F}_j}{\alpha_j}
  \end{aligned}
\end{equation*}
and 
\begin{equation*}
\begin{aligned}
  {\widetilde{u}_{\jph}^R}&=\xbar{u}_\jph^n-\frac{\mathcal{F}_{j+1}-f(\xbar{u}_\jph^n)}{\alpha_{j+1}}
  =\frac{\xbar{u}_{j+\frac{3}{2}}^n+\xbar{u}_\jph^n}{2}-\frac{f(\xbar{u}_{j+\frac{3}{2}}^n)-f(\xbar{u}_\jph^n)}{2\alpha_{j+1}}-{\eta_{j+1}\frac{\Delta\mathcal{F}_{j+1}}{\alpha_{j+1}}}\\
  &=u_{j+1}^*-\eta_{j+1}\frac{\Delta\mathcal{F}_{j+1}}{\alpha_{j+1}},
  \end{aligned}
\end{equation*}
where $u_j^*$ and $u_{j+1}^*$ are nothing but the first-order finite volume Riemann intermediate states, which takes the following form 
\begin{equation}\label{u_star1}
    u_j^*=\frac{\xbar{u}_\jmh^n+\xbar{u}_\jph^n}{2}-\frac{f(\xbar{u}_\jph^n)-f(\xbar{u}_\jmh^n)}{2\alpha_j},\quad \forall j.
\end{equation}
As one can observe from \eqref{u_ave_new}, $\xbar{u}_\jph^{n+1}$ is indeed a convex combination of quantities from the previous time step under the standard CFL condition 
\begin{equation*}
    \dt\leq\frac{\dx}{\alpha_j+\alpha_{j+1}},
\end{equation*}
which was also used in \cite{duan2024AF} for the BP property of cell averages and was inspired by \cite{Guermond_IDP2,Vilar_DGFV}. Moreover, since 
\begin{equation*}
    \alpha_j=\alpha(\xbar{u}_\jmh^n,\xbar{u}_\jph^n)=\max\big(\vert f'(\xbar{u}_\jmh^n)\vert,\vert f'(\xbar{u}_\jph^n)\vert\big),
\end{equation*}
it follows that
\begin{equation*}
    u_j^*\in \big[\min(\xbar{u}_\jmh^n,\xbar{u}_\jph^n),\max(\xbar{u}_\jmh^n,\xbar{u}_\jph^n)\big].
\end{equation*}

We now introduce the definition of the blending coefficient to ensure that if the numerical initial solution lies in $\mathcal{D}$, the solution $\xbar{u}_\jph^n$ remains in $\mathcal{D}$ throughout the entire calculation. 
{
\begin{proposition}
In the numerical scheme \eqref{u_ave_new}, if we set 
\begin{equation}\label{eq:theta_average}
  \eta_j=\min\Big(1,\frac{\alpha_j}{\big\vert\Delta \mathcal{F}_j\big\vert}\min\big(\mathring u_{\max}-u_j^*,u_j^*-\mathring{u}_{\min}\big)\Big),\quad \forall j,
\end{equation}
the scheme preserves the following property: if $\xbar u_\jph^n\in\mathcal D$, then $\xbar u_\jph^{n+1}\in\mathcal{D}$, for any $n$ and any $j$.
\end{proposition}}

\begin{proof}
   We need to determine under which condition, {$u_j^*\pm\eta_j\frac{\Delta\mathcal{F}_j}{\alpha_j}\in [\mathring{u}_{\min}, \mathring u_{\max}]$}, knowing that if $\xbar{u}_\jph^n\in [{\mathring{u}_{\min}, \mathring u_{\max}}]$, for all $j$.
   This requires
   \begin{equation}\label{aveIDP}
       {\mathring{u}_{\min}}\leq u_{j}^*+\eta_{j}\frac{\Delta\mathcal{F}_{j}}{\alpha_{j}}\leq {\mathring u_{\max}},\quad {\mathring{u}_{\min}}\leq u_{j}^*-\eta_{j}\frac{\Delta\mathcal{F}_{j}}{\alpha_{j}}\leq {\mathring u_{\max}},\quad \forall j
   \end{equation}
    i.e.,
    \begin{equation*}
        \alpha_{j}\big ({\mathring{u}_{\min}}-u_{j}^*\big ) \leq \eta_{j}\,\Delta\mathcal{F}_{j}\leq \alpha_{j}\big ({\mathring{u}_{\max}}-u_{j}^*\big ),~ \alpha_{j}\big (u_{j+1}^*-{\mathring{u}_{\min}}\big ) \geq \eta_{j}\,\Delta\mathcal{F}_{j}\geq \alpha_{j}\big (u_{j}^*-{\mathring{u}_{\max}}\big ),~ \forall j.
    \end{equation*}
    Since
    \begin{equation*}
        \alpha_{j}\big ( {\mathring{u}_{\min}}-u_{j}^*\big )\leq 0\leq \alpha_{j}\big ( {\mathring{u}_{\max}}-u_{j}^*\big ),\quad \alpha_{j}\big ( u_{j}^*-{\mathring{u}_{\min}}\big )\geq 0\geq \alpha_{j}\big ( u_{j}^*-{\mathring{u}_{\max}}\big ),\quad \forall j,
    \end{equation*}
    we have two cases for $\eta_j$:
    \begin{itemize}
        \item If $\Delta\mathcal{F}_{j}\geq 0$, then we need
        \begin{equation*}
            \eta_{j}\leq \alpha_{j}\dfrac{\big ({\mathring{u}_{\max}}-u_{j}^*\big ) }{\Delta\mathcal{F}_{j}}\quad\mbox{and}\quad  \eta_{j}\leq \alpha_{j}\dfrac{\big (u_{j}^*-{\mathring{u}_{\min}}\big ) }{\Delta\mathcal{F}_{j}};
        \end{equation*}
        \item If $\Delta\mathcal{F}_{j}\leq 0$, then we need
        \begin{equation*}
            \eta_{j}\leq \alpha_{j}\dfrac{\big (u_{j}^*-{\mathring{u}_{\min}}\big ) }{-\Delta\mathcal{F}_{j}}\quad\mbox{and}\quad  \eta_{j}\leq \alpha_{j}\dfrac{\big ({\mathring{u}_{\max}}-u_{j}^*\big ) }{-\Delta\mathcal{F}_{j}}.
        \end{equation*}
    \end{itemize}
    In addition, we require $\eta_{j}\leq 1$. Therefore, we get the result as given in \eref{eq:theta_average}.
\end{proof}

\subsection{Blended residuals and update of point value}\label{subsec:point}
The point value will be evolved in a similar convex limiting manner as described in the previous subsection. To be specific, each point value is assigned two residuals, one is the low-order residual \eref{residu_LO} and one is the high-order residual \eref{residuals_HO}. These two residuals will be blended in a convex manner through blending coefficients $\theta_\jmh\in[0,1]$ and $\theta_\jph\in[0,1]$. The blended residuals read as  
\begin{equation}\label{eq:blend_residuals}
\left\{
\begin{aligned}
\overrightarrow{\Phi}_\jmh
  &=\overrightarrow{\Phi}_\jmh^{LO}+\theta_\jmh(\overrightarrow{\Phi}_\jmh^{HO}-\overrightarrow{\Phi}_\jmh^{LO})
  =\overrightarrow{\Phi}_\jmh^{LO}+\theta_\jmh\Delta\overrightarrow{\Phi}_\jmh,\\
  \overleftarrow{\Phi}_\jph
  &=\overleftarrow{\Phi}_\jph^{LO}+\theta_\jph(\overleftarrow{\Phi}_\jph^{HO}-\overleftarrow{\Phi}_\jph^{LO})
  =\overleftarrow{\Phi}_\jph^{LO}+\theta_\jph\Delta\overleftarrow{\Phi}_\jph.
  \end{aligned}\right.
\end{equation}
Analogously, when the blending coefficient $\theta_{j\pm\frac{1}{2}}=0$, it results in the first-order residuals while if $\theta_{j\pm\frac{1}{2}}=1$, it induces to the high-order case. Since conservation is not needed for point value, we can assign two different blending coefficients $\theta_{j\pm\frac{1}{2}}$ for the residuals that {are} used for the update of point value at the cell interface $x=x_j$.

The first-order forward Euler time integration for point value is expressed as
\begin{equation}\label{eq:point_new}
  u_j^{n+1}=u_j^n-2\lambda(\overrightarrow{\Phi}_\jmh+\overleftarrow{\Phi}_\jph),\quad \forall j,
\end{equation}
and then it is rewritten as a convex combination of quantities defined at the previous time step, i.e.,
\begin{equation}\label{u_point_new}
\begin{aligned}
u_j^{n+1}&=u_j^n-2\lambda(\overrightarrow{\Phi}_\jmh+\overleftarrow{\Phi}_\jph)+2\lambda(\beta_\jmh+\beta_\jph) u_j^n- 2\lambda(\beta_\jmh+\beta_\jph) u_j^n\\
&=\Big(1-2\lambda\beta_\jmh-2\lambda\beta_\jph\Big)u_j^n
 +2\lambda\beta_\jmh\Big(u_j^n-\frac{\overrightarrow{\Phi}_\jmh}{\beta_\jmh}\Big)
 +2\lambda\beta_\jph\Big(u_j^n-
 \frac{\overleftarrow{\Phi}_\jph}{\beta_\jph}\Big)\\
&=\Big(1-2\lambda\beta_\jmh-2\lambda\beta_\jph\Big)u_j^n +2\lambda\beta_\jmh {\widetilde u_j^L}+2\lambda\beta_\jph{\widetilde u_j^R},
\end{aligned}  
\end{equation}
where
\begin{equation}\label{theta_1}
\begin{aligned}
  {\widetilde u_j^L}&=u_j^n-\frac{\overrightarrow{\Phi}_\jmh}{\beta_\jmh}=\frac{\xbar u_\jmh^n+u_j^n}{2}-\frac{f(u_j^n)-f(\xbar u_\jmh^n)}{2\beta_\jmh}-\theta_\jmh\frac{\Delta \overrightarrow{\Phi}_\jmh}{\beta_\jmh}\\
  &={u^{*,L}_j}-\theta_\jmh\frac{\Delta \overrightarrow{\Phi}_\jmh}{\beta_\jmh}
\end{aligned}
\end{equation}
and
\begin{equation}\label{theta_2}
\begin{aligned}
 {\widetilde u_j^R}&=u_j^n-
\frac{\overleftarrow{\Phi}_\jph}{\beta_\jph}=\frac{u_j^n+\xbar u_\jph^n}{2}-\frac{f(\xbar u_\jph^n)-f(u_j^n)}{2\beta_\jph}-\theta_\jph\frac{\Delta\overleftarrow{\Phi}_\jph}{\beta_\jph}\\
&={u_j^{*,R}}-\theta_\jph\frac{\Delta\overleftarrow{\Phi}_\jph}{\beta_\jph}.
\end{aligned}
\end{equation}
In the above, we have defined 
\begin{equation}\label{u_star2}
  u_j^{*,L}=\frac{\xbar u_\jmh^n+u_j^n}{2}-\frac{f(u_j^n)-f(\xbar u_\jmh^n)}{2\beta_\jmh},\quad u_j^{*,R}=\frac{u_j^n+\xbar u_\jph^n}{2}-\frac{f(\xbar u_\jph^n)-f(u_j^n)}{2\beta_\jph},\quad \forall j.
\end{equation}

We then perform the similar analysis as in subsection \ref{subsec:average} to satisfy the condition of convex combination and obtain:
{
\begin{proposition}
  If we set 
\begin{equation}\label{eq:theta:point}
\theta_\jmh=
  \min\Big(1,\frac{\beta_\jmh}{\big\vert\Delta \overrightarrow{\Phi}_\jmh\big\vert}\min\big(u_j^{*,L}-\mathring{u}_{\min},\mathring{u}_{\max}-u_j^{*,L}\big)\Big)
\end{equation}
and
\begin{equation}\label{eq:theta:point2}
\theta_\jph=
\min\Big(1,\frac{\beta_\jph}{\big\vert\Delta \overleftarrow{\Phi}_\jph\big\vert}\min\big(u_j^{*,R}-\mathring{u}_{\min},\mathring{u}_{\max}-u_j^{*,R}\big)\Big),
\end{equation}
the numerical scheme  \eqref{u_point_new} has the same BP property: if $u_j^n\in\mathcal D$, then $u_j^{n+1}\in\mathcal{D}$, for any $n$ and any $j$.
\end{proposition}}

\begin{proof}
The proof is very similar to that for \eqref{eq:theta_average}. We assume that
\begin{equation}\label{pointIDP}
    {\mathring{u}_{\min}}\leq u_j^{*,L} -\theta_\jmh\frac{\Delta\overrightarrow{\Phi}_\jmh}{\beta_\jmh}\leq {\mathring{u}_{\max}}, \quad  {\mathring{u}_{\min}}\leq u_j^{*,R} -\theta_\jph\frac{\Delta\overleftarrow{\Phi}_\jph}{\beta_\jph}\leq {\mathring{u}_{\max}},
\end{equation}
and since
\begin{equation*}
  \beta_\jmh \big ( u_j^{*,L}-{\mathring{u}_{\max}}\big )\leq 0\leq \beta_\jmh\big (u_j^{*,L}-{\mathring{u}_{\min}}\big ),~\beta_\jph \big ( u_j^{*,R}-{\mathring{u}_{\max}}\big )\leq 0\leq \beta_\jph\big (u_j^{*,R}-{\mathring{u}_{\min}}\big ),  
\end{equation*}
we can rewrite \eref{pointIDP} as
\begin{equation*}
\begin{aligned}
    &\beta_\jmh\big (u_j^{*,L}-{\mathring{u}_{\max}}\big )\leq \theta_\jmh\Delta\overrightarrow{\Phi}_\jmh\leq \beta_\jmh \big ( u_j^{*,L}-{\mathring{u}_{\min}}\big ), \\ &\beta_\jph\big (u_j^{*,R}-{\mathring{u}_{\max}}\big )\leq \theta_\jph\Delta\overleftarrow{\Phi}_\jph\leq \beta_\jph \big ( u_j^{*,R}-{\mathring{u}_{\min}}\big ).
    \end{aligned}
\end{equation*}
We then have two cases for $\theta_\jmh$:
\begin{itemize}
        \item If $\Delta\overrightarrow{\Phi}_\jmh\geq 0$, then we need
        \begin{equation*}
          \theta_\jmh\leq \dfrac{\beta_\jmh\big (u_j^{*,L}-{\mathring{u}_{\min}}\big )}{\Delta\overrightarrow{\Phi}_\jmh};  
        \end{equation*}
        \item If $\Delta\overrightarrow{\Phi}_\jmh\leq 0$, then we need
        \begin{equation*}
           \theta_\jmh\leq \dfrac{\beta_\jmh\big ({\mathring{u}_{\max}}-u_j^{*,L}\big )}{-\Delta\overrightarrow{\Phi}_\jmh}. 
        \end{equation*}
    \end{itemize}
    In case of an inaccurate estimation of the upwind direction which may lead to a opposite sign of $\Delta\overrightarrow{\Phi}_\jmh$, especially in the transonic case and $\theta_\jmh\leq1$ is also required, we therefore take
    \begin{equation*}
       \theta_\jmh=
  \min\Big(1,\frac{\beta_\jmh}{\big\vert\Delta \overrightarrow{\Phi}_\jmh\big\vert}\min\big(u_j^{*,L}-{\mathring{u}_{\min}},{\mathring{u}_{\max}}-u_j^{*,L}\big)\Big). 
    \end{equation*}
    Similar analysis can be performed for $\theta_\jph$, we omit here for saving space.
\end{proof}

\subsection{Local maximum principle and smooth extrema relaxation}\label{subsec33} 
In \eref{eq:theta_average}, \eref{eq:theta:point}, and \eref{eq:theta:point2}, the global bounds are taken as 
\begin{equation}\label{GLMP}
    {\mathring{u}_{\min}}=\min_{x\in\Omega}\; u_0(x) \quad \mbox{and} \quad {\mathring{u}_{\max}}=\max_{x\in\Omega}\; u_0(x).
\end{equation}
In order to reduce as much as possible the appearance of spurious oscillations in the approximation of discontinuous solution, in addition to the global maximum principle, we also impose a local maximum principle for the update of average value. More specific, we take the local bounds of the solution from the previous time step (or the previous Runge-Kutta step in the general case) and we set
\begin{equation}\label{Local_MP}
    u_\jph^{\min}=\min\big(\xbar{u}_\jph^n, u_j^*,u_{j+1}^*\big)\quad \mbox{and}\quad u_\jph^{\max}=\max\big(\xbar{u}_\jph^n, u_j^*,u_{j+1}^*\big),\quad \forall j,
\end{equation}
for the update of cell average. To ensure the updated cell average bounded by these local bounds and in the light of the combination in \eref{u_ave_new}, we need to guarantee both ${\widetilde{u}_\jph^L}\in[u_\jph^{\min},u_\jph^{\max}]$ and ${\widetilde{u}_\jmh^R}\in[u_\jmh^{\min},u_\jmh^{\max}]$, it is sufficient to take $\eta_j$ as follows:
\begin{equation}\label{theta_loc:ave}
\eta_j=\left\{
\begin{aligned}
&\min\Big(1,\frac{\alpha_j}{\vert\Delta \mathcal{F}_j\vert}\min(u_\jph^{\max}-u_j^*,u_j^*-u_\jmh^{\min})\Big)&&\mbox{if}~\Delta \mathcal{F}_j>0,\\
&\min\Big(1,\frac{\alpha_j}{\vert\Delta \mathcal{F}_j\vert}\min(u_j^*-u_\jph^{\min},u_\jmh^{\max}-u_j^*)\Big)&&\mbox{if}~\Delta \mathcal{F}_j<0.
\end{aligned}\right.
\end{equation}

\begin{proof}
 We need to determine under which condition, $u_j^*+\eta_j\frac{\Delta\mathcal{F}_j}{\alpha_j}\in [u_\jph^{\min}, u_\jph^{\max}]$ and $u_{j}^*-\eta_{j}\frac{\Delta\mathcal{F}_{j}}{\alpha_{j}}\in [u_\jmh^{\min}, u_\jmh^{\max}]$, for all $j$. Since $u_j^*\in[u_\jph^{\min}, u_\jph^{\max}]$ and $u_j^*\in[u_\jmh^{\min}, u_\jmh^{\max}]$, we have two cases for $\eta_j$:
    \begin{itemize}
        \item If $\Delta\mathcal{F}_{j}\geq 0$, then we need $u_j^*+\eta_j\frac{\Delta\mathcal{F}_j}{\alpha_j}\leq u_\jph^{\max}$ and $u_j^*-\eta_j\frac{\Delta\mathcal{F}_j}{\alpha_j}\geq u_\jmh^{\min}$, which yields
        \begin{equation*}
            \eta_{j}\leq \alpha_{j}\dfrac{\big (u_\jph^{\max}-u_{j}^*\big ) }{\Delta\mathcal{F}_{j}}\quad\mbox{and}\quad  \eta_{j}\leq \alpha_{j}\dfrac{\big (u_{j}^*-u_\jmh^{\min}\big ) }{\Delta\mathcal{F}_{j}};
        \end{equation*}
        \item If $\Delta\mathcal{F}_{j}\leq 0$, then we need $u_j^*+\eta_j\frac{\Delta\mathcal{F}_j}{\alpha_j}\geq u_\jph^{\min}$ and $u_j^*-\eta_j\frac{\Delta\mathcal{F}_j}{\alpha_j}\leq u_\jmh^{\max}$, which yields
        \begin{equation*}
            \eta_{j}\leq \alpha_{j}\dfrac{\big (u_{j}^*-u_\jph^{\min}\big ) }{-\Delta\mathcal{F}_{j}}\quad\mbox{and}\quad  \eta_{j}\leq \alpha_{j}\dfrac{\big (u_\jmh^{\max}-u_{j}^*\big ) }{-\Delta\mathcal{F}_{j}}.
        \end{equation*}
    \end{itemize}
    In addition, we require $\eta_{j}\leq 1$. Therefore, we get the result as given in \eref{theta_loc:ave}. 
\end{proof}

{The formulation \eqref{theta_loc:ave} for the blending parameter is standard (see, e.g., \cite{Kuzmin_MCL1,Vilar_DGFV}) and was also used in \cite{duan2024AF}.}

\paragraph*{Smooth extrema relaxation}
In order to preserve high-accuracy in the vicinity of smooth extrema, we employ the smooth detector used in \cite{AbgrallLiu, Vilar} to relax the local maximum principle constraints for the average value update. This detector is designed as follows: 
         \begin{itemize}
        \item First, we denote the parabolic approximation in \eref{u_para} by $\mathfrak{u}_\jph$.
           \item Next, we compute $\zeta'=\mathfrak{u}_\jph'(x_\jph)$, $\zeta'_L=\mathfrak{u}_\jph'(x_j)$, and the minimum and maximum values of the derivative around cell $I_\jph$:
               \begin{equation*}
                 \zeta'^{,\jph}_{\min/\max}=\min/\max(\mathfrak{u}'_\jmh(x_i),
                 \mathfrak{u}'_{j+\frac{3}{2}}(x_i)),\quad i=j, j+1.
               \end{equation*}
           All of the obtained quantities are then used to define the left detection factor $\xi_L$:
           \begin{equation*}
              \xi_L=\left\{\begin{aligned}
              &\min\Big(1,\frac{\zeta'^{,\jph}_{\max}-\zeta'}{\zeta'_L-\zeta'}\Big),\quad&& \mbox{if}~\zeta'_L>\zeta',\\
              &1,\quad && \mbox{if}~\zeta'_L=\zeta',\\
              &\min\Big(1,\frac{\zeta'^{,\jph}_{\min}-\zeta'}{\zeta'_L-\zeta'}\Big),\quad&& \mbox{if}~\zeta'_L<\zeta'.
             \end{aligned}\right.
            \end{equation*}
           \item Analogously, we compute $\zeta'=\mathfrak{u}_\jph'(x_\jph)$, $\zeta'_R=\mathfrak{u}_\jph'(x_{j+1})$ and define the right detection factor $\xi_R$ as
                \begin{equation*}
              \xi_R=\left\{\begin{aligned}
              &\min\Big(1,\frac{\zeta'^{,\jph}_{\max}-\zeta'}{\zeta'_R-\zeta'}\Big),\quad&& \mbox{if}~\zeta'_R>\zeta',\\
              &1,\quad && \mbox{if}~\zeta'_R=\zeta',\\
             &\min\Big(1,\frac{\zeta'^{,\jph}_{\min}-\zeta'}{\zeta'_R-\zeta'}\Big),\quad&& \mbox{if}~\zeta'_R<\zeta'.
              \end{aligned}\right.
           \end{equation*}
           \item Finally, if $\xi_L=1$ and $\xi_R=1$, we consider that the numerical solution presents a smooth profile on cell $[x_j,x_{j+1}]$. In this particular case, the blended coefficient constraint for the average value update through the local maximum principle is relaxed.
        \end{itemize}
        
The basic idea behind this detector is as follows: the numerical solution is expected to exhibit a smooth extremum if at least the linearized spatial derivatives of the numerical solution present a monotonic profile. We aim to preserve scheme accuracy in the presence of such a smooth extremum and verify whether the gradient at the cell boundaries lies within the interval $[\zeta'^{,\jph}_{\min}, \zeta'^{,\jph}_{\max}]$.

In the end, we would like to remark that we neither apply any local maximum principles nor design any smooth extrema detectors for the update of point value, since we have less requirements/restrictions on the point value and from the numerical experiments conducted in section \ref{sec4}, we do not observe any issues caused by this. The smooth extrema detector for update of the point value will need further investigation.

\bigskip 
{
We end this section by presenting in Algorithm \ref{alg1} a pseudo-code for the forward Euler evolution of the proposed BP PAMPA scheme for scalar conservation laws.
\begin{algorithm}
\caption{{Forward Euler evolution of BP PAMPA scheme}}
\label{alg1}
\renewcommand{\algorithmicrequire}{\textbf{Input:}}
\renewcommand{\algorithmicensure}{\textbf{Output:}}
\begin{algorithmic}[1]
\REQUIRE{$\xbar{u}_\jph$, $u_j$ at time $t=t^n$.}
\STATE{Compute the parameters $\alpha_j$ in \eqref{alpha_ave} and $\beta_{j\pm\frac{1}{2}}$ in \eqref{beta_point}.}
\STATE{Compute the flux blending coefficient in \eqref{eq:theta_average} using \eqref{u_star1} and \eqref{LLF_flux}.}
\STATE{Check the smooth extrema relaxation in subsection \ref{subsec33} }
\IF {$\xi_L\neq1$ or $\xi_R\neq1$}
\STATE{Apply the local maximum principle and recompute the flux blending coefficient using \eqref{Local_MP} and \eqref{theta_loc:ave}.}
\ENDIF
\STATE{Compute the residuals blending coefficients in \eqref{eq:theta:point} and \eqref{eq:theta:point2}, using \eqref{residu_LO}, \eqref{residuals_HO}, and \eqref{u_star2}.}
\STATE{Compute the fluctuations. Namely,
\begin{itemize}
  \item Compute blended fluxes using \eqref{eq:blend_flux}.
  \item Compute blended residuals using \eqref{eq:blend_residuals}.
\end{itemize}
}
\STATE{Evolve $\xbar u_\jph$ and $u_j$ simultaneously using \eqref{eq:ave_new} and \eqref{eq:point_new}, respectively.}
\ENSURE{$\xbar u_\jph$ and $u_j$ at the new time level $t^{n+1}$.}
\end{algorithmic}
\end{algorithm}}

\subsection{Case of the Euler equations}\label{EulerGQL} 
For the Euler equations of gas dynamics, the invariant domain $\mathcal{D}$ is given by \eref{Euler_IDP}. Following \cite{wu2023geometric}, we see that the invariant domain is exactly equivalent to the set given in \eref{Euler_GQL}, which can be rewritten as
\begin{equation}\label{GQL_D}
    \mathcal{D}_\nu=\{\mathbf u=(\rho, m, E)^T: \rho> 0 ~\text{and for all}~\nu\in \mathbb{R}, \mathbf u\cdot \bm \nu>0\}, ~\bm \nu=(\frac{\nu^2}{2},-\nu,1)^T.
\end{equation}
The first relation is obvious, while the second one comes from:
\begin{equation*}
   \rho\frac{\nu^2}{2}-\nu m+E=e+\frac{\rho}{2}(u-\nu)^2, 
\end{equation*}
where $e=E-\tfrac{m^2}{2\rho}$ and $m=\rho v$. Note that the constraints in the GQL representation \eref{GQL_D} become linear with respect to $\mathbf u$. The goal is reduced to find optimal blending coefficients to satisfy the linear constraints. 

\paragraph*{Positivity preserving} We now illustrate the preservation of positive density and internal energy. Using the same notations as in the scalar case (see \eref{aveIDP} and \eref{pointIDP}), the scheme is positivity preserving if the following constraints are fulfilled:
\begin{itemize}
\item $\rho_j^*+\eta_j^{\rho}\dfrac{\Delta\mathcal{F}_j^{\rho}}{\alpha_j}> 0$ and $\rho_j^*-\eta_j^{\rho}\dfrac{\Delta\mathcal{F}_j^{\rho}}{\alpha_j}> 0$;
\item $\rho_\jmh^*-\theta_\jmh^{\rho}\dfrac{\Delta\overrightarrow{\Phi}_\jmh^{\rho}}{\beta_\jmh}> 0$;
\item $\rho_\jph^*-\theta_\jph^{\rho}\dfrac{\Delta\overleftarrow{\Phi}_\jph^{\rho}}{\beta_\jph}> 0$;
\item $\mathbf u_{j}^*\cdot\bm\nu+\eta_j^e\dfrac{\Delta\mathcal{F}_{j}\cdot\bm\nu}{\alpha_{j}}> 0$ and $\mathbf u_{j}^*\cdot\bm\nu-\eta_j^e\dfrac{\Delta\mathcal{F}_{j}\cdot\bm\nu}{\alpha_{j}}> 0$;
\item $\mathbf u_\jmh^*\cdot\bm \nu-\theta_\jmh^e\dfrac{\Delta \overrightarrow{\Phi}_\jmh\cdot \bm\nu}{\beta_\jmh}> 0$;
\item $\mathbf u_\jph^*\cdot \bm\nu-\theta_\jph^e\dfrac{\Delta \overrightarrow{\Phi}_\jph\cdot\bm\nu}{\beta_\jph}> 0$.
\end{itemize}

For the density, according to \eref{eq:theta_average}, we obtain 
\begin{equation*}
    \eta_j^{\rho}=\min\Big(1,\frac{\alpha_j}{\big\vert\Delta \mathcal{F}_j^\rho\big\vert}\rho_j^*\Big),
\end{equation*}
while  according to \eref{eq:theta:point} and \eref{eq:theta:point2}, we have
\begin{equation*}
  \theta_\jmh^{\rho}=\min\Big(1,\frac{\beta_\jmh}{\big\vert\Delta \overrightarrow{\Phi}_\jmh^\rho\big\vert}\rho_\jmh^* \Big)\quad\mbox{and}\quad  \theta_\jph^{\rho}=
\min\Big(1,\frac{\beta_\jph}{\big\vert\Delta \overleftarrow{\Phi}_\jph^\rho\big\vert}\rho_\jph^*\Big).
\end{equation*}

For the internal energy, we take \eref{eq:theta_average} with
$\Delta\mathcal{F}_j$ replaced by $\Delta\mathcal{F}_j\cdot\bm \nu$, while take \eref{eq:theta:point} and \eref{eq:theta:point2} with $\Delta \overrightarrow{\Phi}_\jmh$ replaced by $\Delta \overrightarrow{\Phi}_\jmh\cdot \bm\nu$ and $\Delta \overleftarrow{\Phi}_\jph$ replaced by $\Delta \overleftarrow{\Phi}_\jph\cdot \bm\nu$, respectively. Then, we need to minimize the expressions that will now depend on $\nu\in \R$. That is, for blending the numerical fluxes, we need to minimize
\begin{equation*}
    \eta_j^e=\alpha_j\min_{\bm\nu}\frac{\mathbf u_j^*\cdot\bm\nu}{\vert\Delta\mathcal{F}_j\cdot\bm\nu\vert},\quad \forall j,
\end{equation*}
and for blending the residuals, we need to minimize
\begin{equation*}
    \theta^e_\jmh=\beta_\jmh\min_{\bm\nu}\frac{\mathbf u_\jmh^*\cdot\bm\nu}{\vert\Delta\overrightarrow{\Phi}_\jmh\cdot\bm\nu\vert},\quad \theta^e_\jph=\beta_\jph\min_{\bm\nu}\frac{\mathbf u_\jph^*\cdot\bm\nu}{\vert\Delta\overleftarrow{\Phi}_\jph\cdot\bm\nu\vert}.
    \end{equation*} 
Any of the scalar product has the form
\begin{equation*}
    P(\nu)=\alpha \frac{\nu^2}{2}-\beta \nu+\gamma,
\end{equation*}
and we need to minimize expressions of the form
\begin{equation*}
    \dfrac{ \mathbf u\cdot \bm{\nu}}{\vert P(\nu)\vert}.
\end{equation*}
This is an expression of the form:
\begin{equation*}
    \varphi(\nu)=\frac{a_1 \frac{\nu^2}{2}-b_1 \nu+c_1}{\vert a_2 \frac{\nu^2}{2}-b_2 \nu+c_2\vert}=
\dfrac{a_1 \nu^2-2b_1 \nu+2c_1}{\vert a_2 \nu^2-2b_2 \nu+2c_2\vert}:=\frac{Q(\nu)}{\vert P(\nu)\vert},
\end{equation*}
where $P(\nu)=a_2 \nu^2-2b_2 \nu+2c_2$ and $Q(\nu)=a_1 \nu^2-2b_1 \nu+2c_1>0$.
Note that in this problem, the parameters are such that $Q(\nu)> 0$ for any $\nu\in\mathbb{R}$. 

In the appendix \ref{appendix:A}, we provide analytical expressions of the minimization problem. These expressions have been used for the numerical test cases. Once the minimization problem is solved, we obtain the blending coefficients, denoted by $\eta_j^e$ or $\theta_{j\pm\frac{1}{2}}^e$, for ensuring the bounds of internal energy. Finally, the blending coefficients used for the evaluation of BP fluxes and residuals are defined as the minimum of the two obtained coefficients, i.e.,
\begin{equation*}
    \eta_j=\min(\eta_j^\rho,\eta_j^e),\quad \theta_{j\pm\frac{1}{2}}=\min(\theta_{j\pm\frac{1}{2}}^\rho, \theta_{j\pm\frac{1}{2}}^e).
\end{equation*}

\paragraph*{Local maximum principle}
For the non-linear system case, we naively apply the local maximum principle introduced for the scalar case in subsection \ref{subsec33} to the conserved density variable and the local bounds are also obtained based on this selected conserved variable. This local maximum principle is then relaxed using the same smooth extrema detector previously applied for the scalar case but now on the conserved density variable. 

\section{Numerical Examples}\label{sec4}

In this section, we present a series of well-known and challenging numerical experiments, ranging from scalar problems to the Euler system of gas dynamics, to demonstrate the robustness and effectiveness of the proposed PAMPA scheme in BP and shock-capturing performance. The time integration is performed using the third-order strong stability-preserving Runge--Kutta (SSP-RK3) method, with a CFL number of $0.2$.

\subsection{Scalar conservation laws}\label{subsec41}
The PDE studied in this subsection is $u_t + f(u)_x = 0$, with the following three flux functions considered:
\begin{itemize}
    \item Linear advection equation: $f(u) = u$;
    \item Burger's equation: $f(u) = \frac{u^2}{2}$;
    \item Buckley--Leverett non-convex case: $f(u) = \frac{4u^2}{4u^2 + (1 - u)^2}$.
\end{itemize}
For the scalar case, we consider both the global maximum principle (GMP) and the relaxed local maximum principle (relaxed LMP). When setting the global bounds, two cases are taken in account. One is the restrict GMP as in \eref{GLMP} and another one is the relaxed global bounds which allows some undershoots and overshoots with magnitude of {$\max\big(10^{-4},10^{-3}(\mathring{u}_{\max}-\mathring{u}_{\min})\big)$ as usually taken in the papers using MOOD approach (see, e.g., \cite{AbgrallLiu,CDL,Loubere_MOOD,Dumbser_MOOD,Dumbser_MOOD2})}.

\subsubsection*{Example 1---Linear advection equation}
In the first example, we consider the linear advection equation. We first test the accuracy of the BP PAMPA method. To this end, we consider the initial data $u_0(x)=\cos(2\pi x)$ on $[0,1]$ with periodic boundary conditions, and compute the solution until a final time of $t=3$ (3 periods). The results are reported in Table \ref{scala_tab1}. It can be seen that the accuracy of the BP-PAMAPA method is somehow degraded to the second-order in the $L^2$ and $L^\infty$ norms, and much less in the $L^1$ norm. {The problem is smooth and the local maximum principle guaranteed parameter computed in \eqref{theta_loc:ave} will not be activated due to the designed smooth extrema detector. Therefore,} the accuracy degraded phenomenon is because of the strict selection of the global bounds {($\mathring{u}_{\min}$ and $\mathring{u}_{\max}$)} in \eref{GLMP}. To verify this, we then apply a relaxed global bound, allowing undershoots and overshoots with a magnitude of {$\max\big(10^{-4},10^{-3}(\mathring{u}_{\max}-\mathring{u}_{\min})\big)$}. The simulation is rerun under this relaxed GMP, and the results are reported in Table \ref{scala_tab2}. As expected, the third-order accuracy is well preserved. 

\begin{table}[ht!]
\caption{\label{scala_tab1} Example 1 (linear advection): Errors and convergence rates obtained by BP PAMAPA using strict GMP.}
\begin{center}
\begin{tabular}{|c||c|c|c|c|c|c|c|}
\hline
&$\dx$& $L^1$-error &rate&  $L^2$-error& rate&$L^\infty$-error& rate \\
\hline
\multicolumn{1}{|c||}{\multirow{4}{*}{\makecell{cell\\ average}}}&$2\;10^{-2}$     &  $5.991\;10^{-4}$&- &  $8.611\;10^{-4}$& -     &        $2.248\;10^{-3}$  & - \\
&$1\;10^{-2}$ &  $1.059\;10^{-4}$&$2.50$&  $1.874\;10^{-4}$& $2.20$&       $5.926\;10^{-4}$  & $1.92$\\
&$5\;10^{-3}$ &  $1.769\;10^{-5}$& $2.58$&  $4.132\;10^{-5}$& $2.18$&       $1.643\;10^{-4}$  & $1.85$\\
&$2.5\;10^{-3}$ &  $2.935\;10^{-6}$& $2.59$&  $9.109\;10^{-6}$& $2.18$&       $4.786\;10^{-5}$  & $1.78$\\
&$1.25\;10^{-3}$ &  $4.818\;10^{-7}$& $2.61$&  $2.004\;10^{-6}$& $2.18$&       $1.396\;10^{-5}$  & $1.78$\\
\hline
\multicolumn{1}{|c||}{\multirow{4}{*}{\makecell{point\\ value}}}&$2\;10^{-2}$     &  $6.678\;10^{-4}$&- &  $9.309\;10^{-4}$& -     &        $2.175\;10^{-3}$  & - \\
&$1\;10^{-2}$ &  $1.135\;10^{-4}$&$2.56$&  $1.981\;10^{-4}$& $2.23$&       $5.983\;10^{-4}$  & $1.86$\\
&$5\;10^{-3}$ &  $1.863\;10^{-5}$& $2.61$&  $4.299\;10^{-5}$& $2.20$&       $1.677\;10^{-4}$  & $1.83$\\
&$2.5\;10^{-3}$ &  $3.041\;10^{-6}$& $2.62$&  $9.366\;10^{-6}$& $2.20$&       $4.887\;10^{-5}$  & $1.78$\\
&$1.25\;10^{-3}$ &  $4.964\;10^{-7}$& $2.61$&  $2.043\;10^{-6}$& $2.20$&       $1.398\;10^{-5}$  & $1.81$\\ \hline
\end{tabular}
\end{center}
\end{table}

\begin{table}[ht!]
\caption{\label{scala_tab2} Same as in Table \ref{scala_tab1} but using the relaxed GMP.}
\begin{center}
\begin{tabular}{|c||c|c|c|c|c|c|c|}
\hline
&$\dx$& $L^1$-error &rate&  $L^2$-error& rate&$L^\infty$-error& rate \\
\hline
\multicolumn{1}{|c||}{\multirow{4}{*}{\makecell{cell\\ average}}}&$2\;10^{-2}$     &  $3.376\;10^{-4}$&- &  $3.959\;10^{-4}$& -     &        $7.848\;10^{-4}$  & - \\
&$1\;10^{-2}$ &  $4.229\;10^{-5}$&$3.00$&  $4.697\;10^{-5}$& $3.08$&       $6.642\;10^{-5}$  & $3.56$\\
&$5\;10^{-3}$ &  $5.290\;10^{-6}$& $3.00$&  $5.875\;10^{-6}$& $3.00$&       $8.309\;10^{-6}$  & $3.00$\\
&$2.5\;10^{-3}$ &  $6.614\;10^{-7}$& $3.00$&  $7.346\;10^{-7}$& $3.00$&       $1.039\;10^{-6}$  & $3.00$\\
&$1.25\;10^{-3}$ &  $8.268\;10^{-8}$& $3.00$&  $9.183\;10^{-8}$& $3.00$&       $1.299\;10^{-7}$  & $3.00$\\
\hline
\multicolumn{1}{|c||}{\multirow{4}{*}{\makecell{point\\ value}}}&$2\,10^{-2}$     &  $3.534\;10^{-4}$&- &  $4.127\;10^{-4}$& -     &        $7.830\;10^{-4}$  & - \\
&$1\;10^{-2}$ &  $4.300\;10^{-5}$&$3.04$&  $4.751\;10^{-5}$& $3.12$&       $6.652\;10^{-5}$  & $3.56$\\
&$5\;10^{-3}$ &  $5.338\;10^{-6}$& $3.01$&  $5.912\;10^{-6}$& $3.01$&       $8.320\;10^{-6}$  & $3.00$\\
&$2.5\;10^{-3}$ &  $6.648\;10^{-7}$& $3.01$&  $7.374\;10^{-7}$& $3.00$&       $1.040\;10^{-6}$  & $3.00$\\
&$1.25\;10^{-3}$ &  $8.295\;10^{-8}$& $3.00$&  $9.207\;10^{-8}$& $3.00$&       $1.300\;10^{-7}$  & $3.00$\\ \hline
\end{tabular}
\end{center}
\end{table}

{
\begin{remark} The blending coefficients in \eqref{eq:theta_average} and \eqref{eq:theta:point}--\eqref{eq:theta:point2} may be strictly less than one, even for smooth problems. This can lead to a degradation of high-order accuracy (see Table \ref{scala_tab1}). The small number $\max\big(10^{-4},10^{-3}(\mathring{u}_{\max}-\mathring{u}_{\min})\big)$ is a parameter used to relax the discrete maximum principle thus allowing for very small undershoots and overshoots, which permits to maintain a good accuracy when dealing with smooth extrema (see Table \ref{scala_tab2}).
\end{remark}
}

\begin{figure}[ht!]
\centerline{\includegraphics[trim=0.8cm 0.4cm 0.9cm 0.5cm,clip,width=5.0cm]{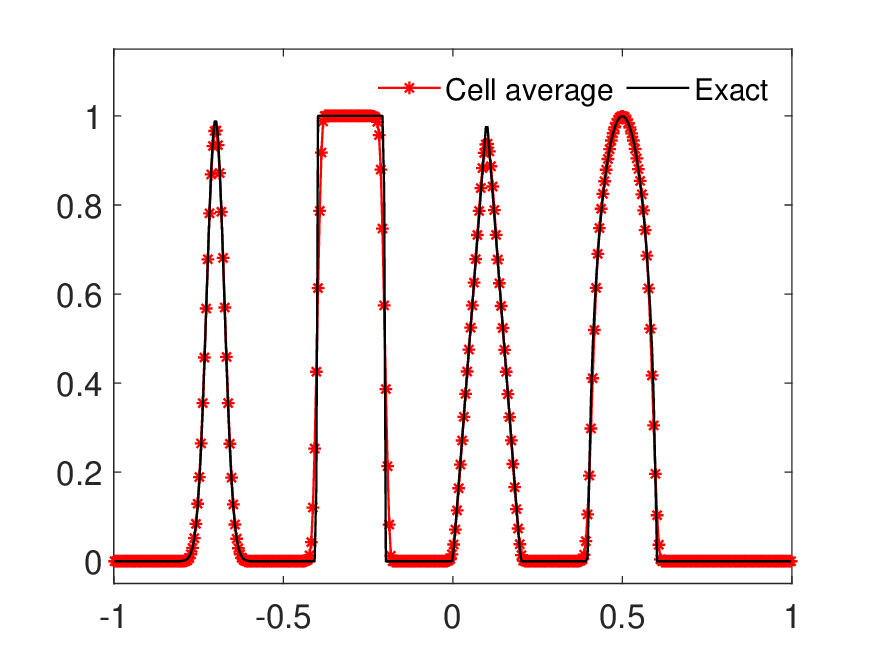}\hspace*{0.5cm}
\includegraphics[trim=0.8cm 0.4cm 0.9cm 0.5cm,clip,width=5.0cm]{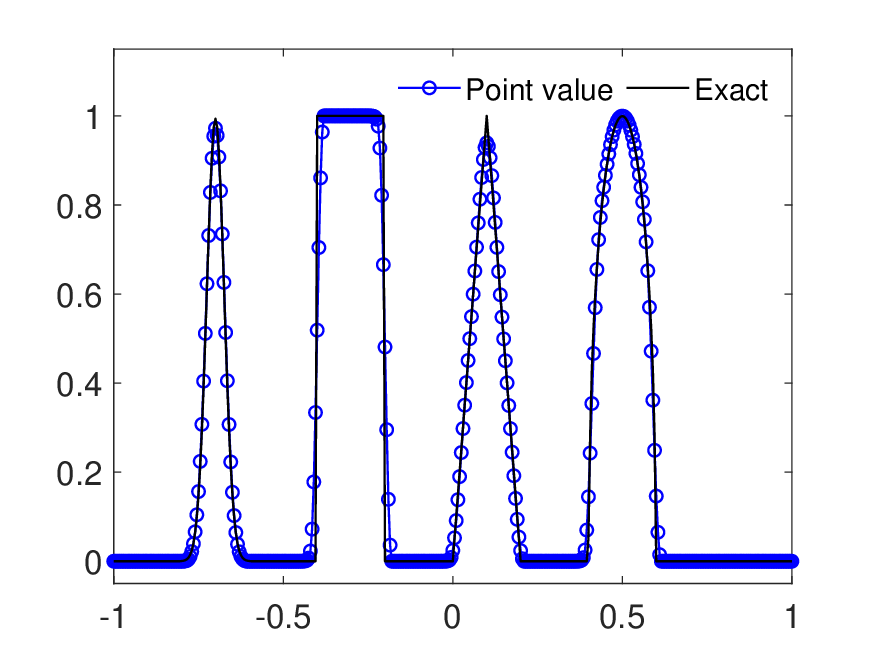}}
\vskip5pt
\centerline{\includegraphics[trim=0.8cm 0.4cm 0.9cm 0.5cm,clip,width=5.0cm]{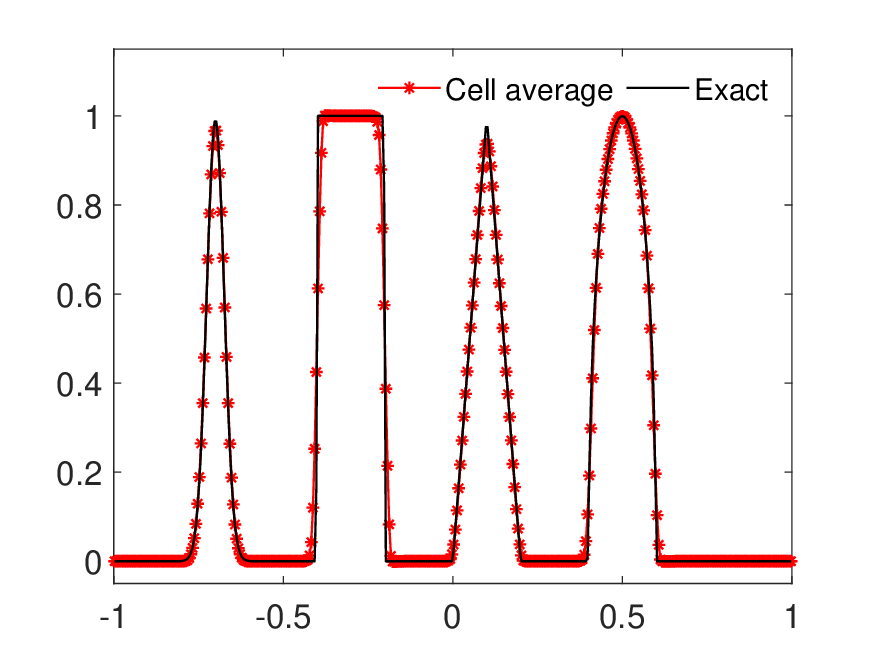}\hspace*{0.5cm}
\includegraphics[trim=0.8cm 0.4cm 0.9cm 0.5cm,clip,width=5.0cm]{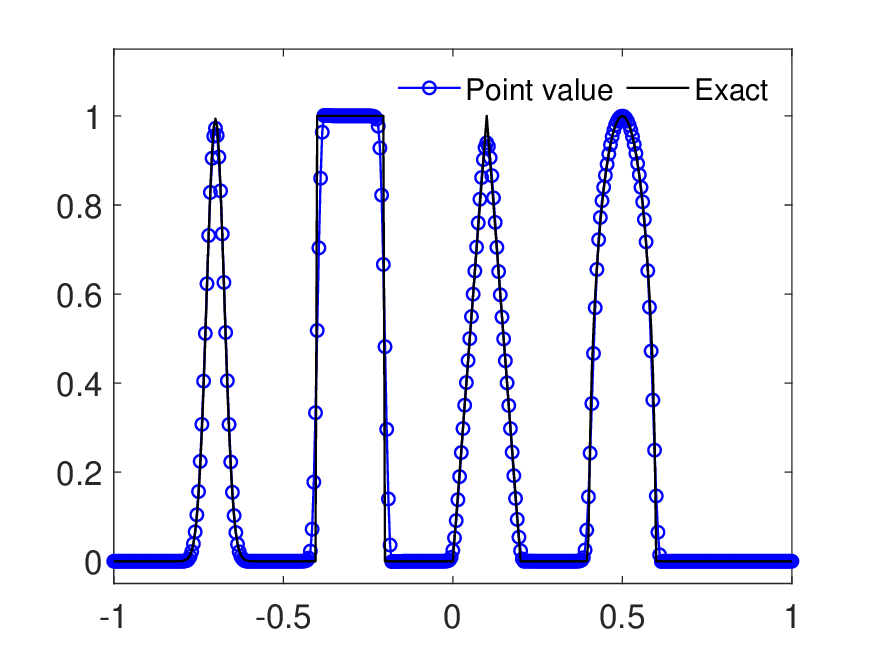}}
\caption{\sf Example 1 (Jiang--Shu): Cell average (left) and point value (right) after one period with $N=400$. Top row: strict GMP $[0,1]$; Bottom row: relaxed GMP $[-0.001, 1.001]$.\label{fig:J-S} }
\end{figure}

Then, we test the non-oscillatory property of the scheme using 
 Jiang--Shu's example \cite{Jiang1996} on the periodic domain $[-1,1]$. The initial data is given by
\begin{equation*}
    u_0(x)=\begin{cases}
    \frac{1}{6}\left(G_1(x, \beta, z-\delta)+G_1(x, \beta, z+\delta)+4 G_1(x, \beta, z)\right) & \text{if}~-0.8 \leqslant x \leqslant-0.6, \\ 
    1 & \text{if}~-0.4 \leqslant x \leqslant-0.2, \\
    1-|10(x-0.1)| & \text{if}~ 0 \leqslant x \leqslant 0.2, \\ \frac{1}{6}\left(G_2(x, \alpha, a-\delta)+G_2(x, \alpha, a+\delta)+4 G_2(x, \alpha, a)\right) & \text{if}~ 0.4 \leqslant x \leqslant 0.6, \\
    0 & \text {else},
    \end{cases}
\end{equation*}
where $G_1(x, \beta, z)=\exp \left(-\beta(x-z)^2\right)$, $G_2(x, \alpha, a)=\sqrt{\max \left(1-\alpha^2(x-a)^2, 0\right)}$. The constants are taken as $a=0.5$, $z=-0.7$, $\delta=0.005$, $\alpha=10$, $\beta=\ln 2 /\left(36 \delta^2\right)$. The solution is composed by the succession of a Gaussian, rectangular, sharp triangular, and parabolic waves. We compute the solution for one period, i.e., until $t=2$, so that the exact solution $u(x,2)=u_0(x)$. The global bounds, used in \eref{eq:theta_average}, \eref{eq:theta:point}, and \eref{eq:theta:point2}, are taken as ${[\mathring{u}_{\min},\mathring{u}_{\max}]}=[0,1]$ for the strict GMP and ${[\mathring{u}_{\min}, \mathring{u}_{\max}]}=[-0.001,1.001]$ for the relaxed GMP. The results computed using $400$ uniform cells are shown in Figure \ref{fig:J-S}, from which one can see how the BP PAMPA scheme behaves, producing a highly accurate solution while ensuring a BP and non-oscillatory profile.

\subsubsection*{Example 2---Non linear Burger's equation}
\begin{figure}[ht!]
\centerline{\includegraphics[trim=0.8cm 0.4cm 0.9cm 0.5cm,clip,width=4.2cm]{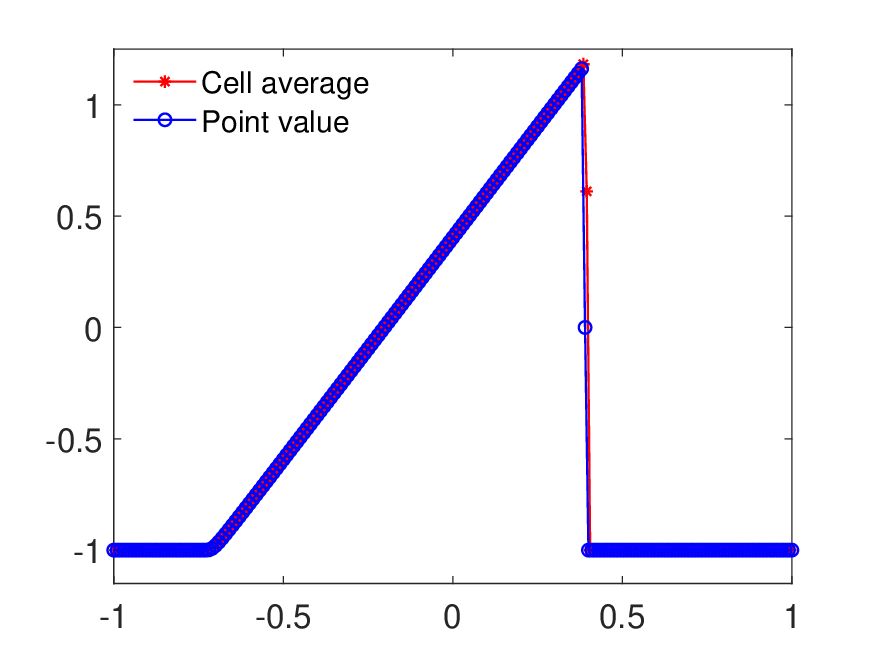}\hspace*{0.05cm}
\includegraphics[trim=0.8cm 0.4cm 0.9cm 0.5cm,clip,width=4.2cm]{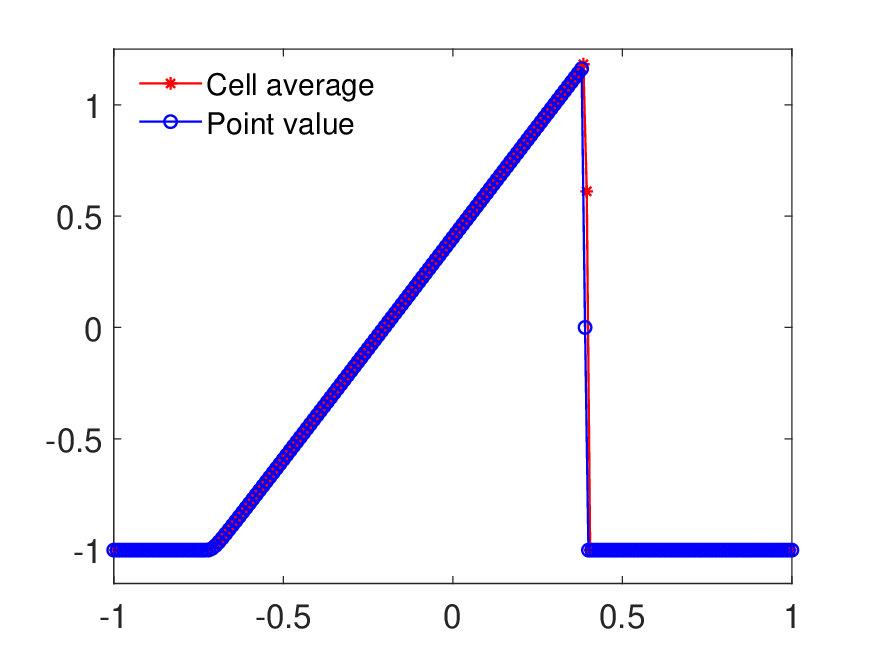}\hspace*{0.05cm}
\includegraphics[trim=0.8cm 0.4cm 0.9cm 0.5cm,clip,width=4.2cm]{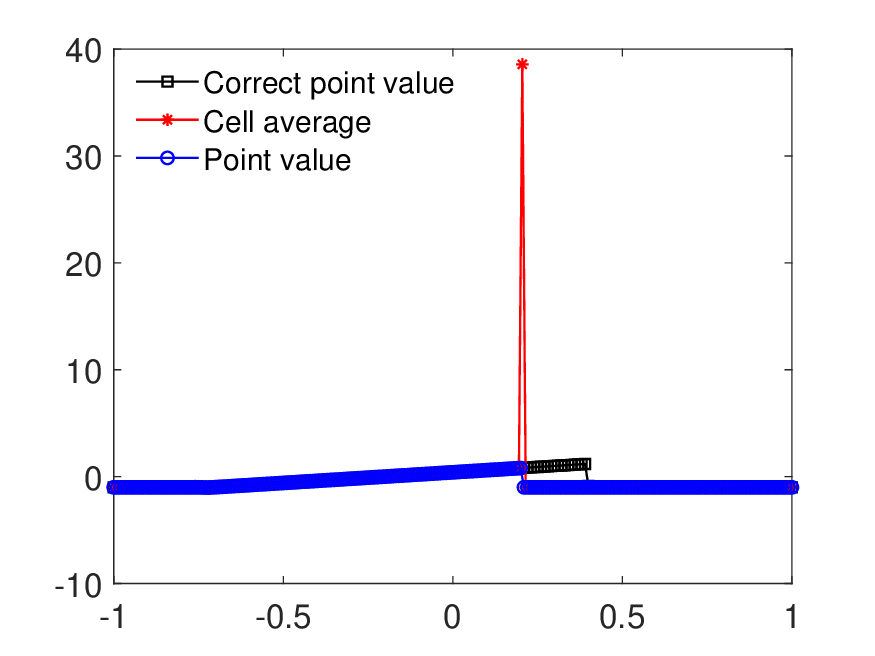}}
\caption{\sf Example 2 (self-steepening shock): Cell average and point value computed by the PAMPA scheme with BP limiting (left and middle, with strict GMP $[-1,2]$ and relaxed GMP $[-1.003,2.003]$, respectively) and without BP limiting (right) {on a uniform mesh with 200 cells}. \label{fig:Burger}}
\end{figure}

In the second example, we consider the nonlinear Burger's equations and test the self-steepening shock on the periodic domain $[-1,1]$ with the initial data given by
\begin{equation*}
    u_0(x)=\left \{ \begin{array}{ll}
2& \text{if}~ \vert x\vert \leq 0.2,\\
-1 & \text{else}.
\end{array}
\right.
\end{equation*}
We compute the solution until the final time $t=0.5$ using $200$ uniform cells. The results obtained by the PAMPA scheme with and without BP limiting are presented in  Figure \ref{fig:Burger}.  A spike is observed in the solution generated by the PAMPA scheme without BP limiting; however, when the proposed BP limiting is applied, no spike appears in the solution, regardless of whether the strict or relaxed GMP is used. The left and middle pictures are consistent with those reported in \cite{HKS,duan2024AF}. The stabilization procedure is enough to avoid the spike shown on the right of Figure \ref{fig:Burger}. It is also worth noting that the left and middle pictures are obtained without sonic point transition that is described in \cite{duan2024AF}.

\subsubsection*{Example 3---Non-convex problem}
\begin{figure}[ht!]
\centerline{\includegraphics[trim=0.8cm 0.4cm 0.9cm 0.5cm,clip,width=5.0cm]{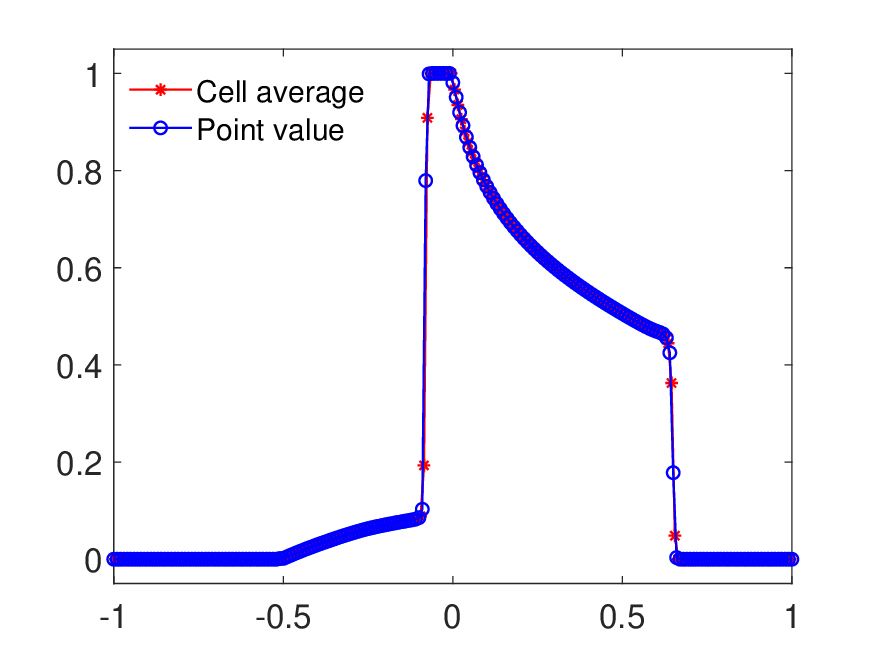}\hspace*{0.5cm}
\includegraphics[trim=0.8cm 0.4cm 0.9cm 0.5cm,clip,width=5.0cm]{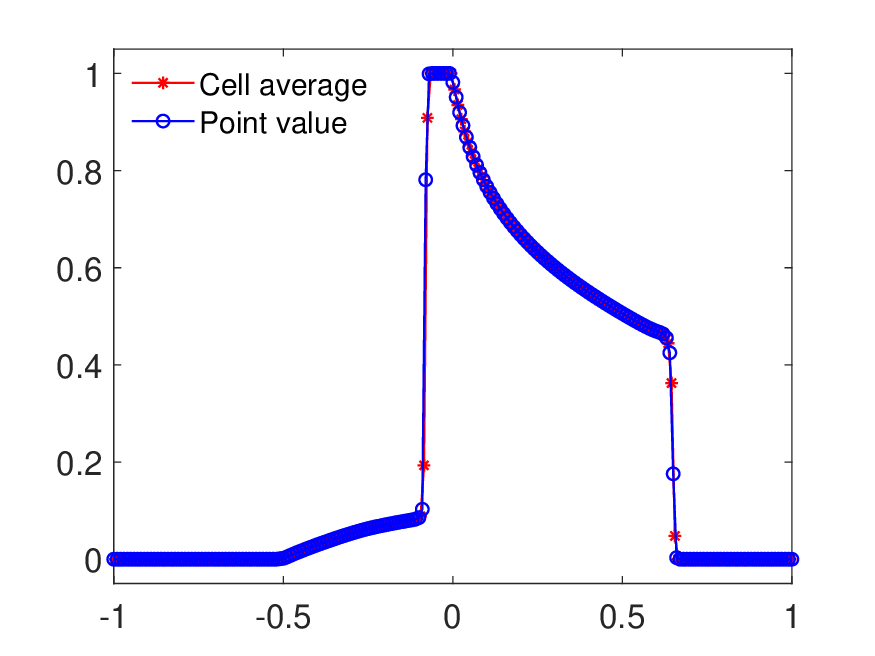}}
\caption{\sf Example 3 (Buckley-Leverett): Cell average and point value computed by the BP PAMPA scheme with restrict GMP (left) and relaxed GMP (right) {on a uniform mesh with 200 cells}. \label{fig:Buckley}}
\end{figure}

In the third example, we test the Buckley--Leverett problem with the initial data defined as
\begin{equation*}
    u_0(x)=\left \{ \begin{array}{ll}
1& \text{if}~ -\frac{1}{2}\leqslant x\leqslant 0,\\
0 & \text{else}.
\end{array}
\right .
\end{equation*}
The computational domain is $[-1,1]$, and the boundary condition is periodic. The numerical solution is computed on a uniform mesh with 200 cells up to the final time $t = 0.4$ and is shown in Figure \ref{fig:Buckley}, demonstrating excellent agreement with the results reported in \cite{Zhang2010}.

\subsection{Compressible Euler equations}
In this subsection, we present a series of challenging numerical examples for the Euler equations of gas dynamics. The specific heat ratio is taken as $\gamma=1.4$, except for Example 4 where $\gamma=3$, and Example 8 where $\gamma=\frac{5}{3}$. For the Euler equations, we consider the preservation of positivity for density and internal energy, as well as the relaxed local maximum principle for the average value of the density variable.

\subsubsection*{Example 4--- Isentropic problem}
In the fourth example, we test a smooth isentropic problem with the initial data given by
\begin{equation*}
    \rho(x,0)=1+0.999995\sin(\pi x),\quad v(x,0)=0,\quad p(x,0)=\rho^{\gamma},\quad x\in[-1,1],
\end{equation*}
and the periodic conditions. 
Because $\gamma=3$, the exact density and velocity in this case can be obtained by the method of characteristics and is explicitly given by
\begin{equation*}
\rho(x,t) = \dfrac12\big( \rho_0(x_1) + \rho_0(x_2)\big), \quad v(x,t) = \sqrt{3}\big(\rho(x,t)-\rho_0(x_1) \big),
\end{equation*}
where for each coordinate $x$ and time $t$ the values $x_1$ and $x_2$ are solutions of the nonlinear equations
\begin{align*}
& x + \sqrt{3}\rho_0(x_1) t - x_1 = 0, \\
& x - \sqrt{3}\rho_0(x_2) t - x_2 = 0.
\end{align*}

\begin{figure}[ht!]
\centerline{\subfigure[cell average of $\rho$]{\includegraphics[trim=0.8cm 0.4cm 0.9cm 0.5cm,clip,width=4.2cm]{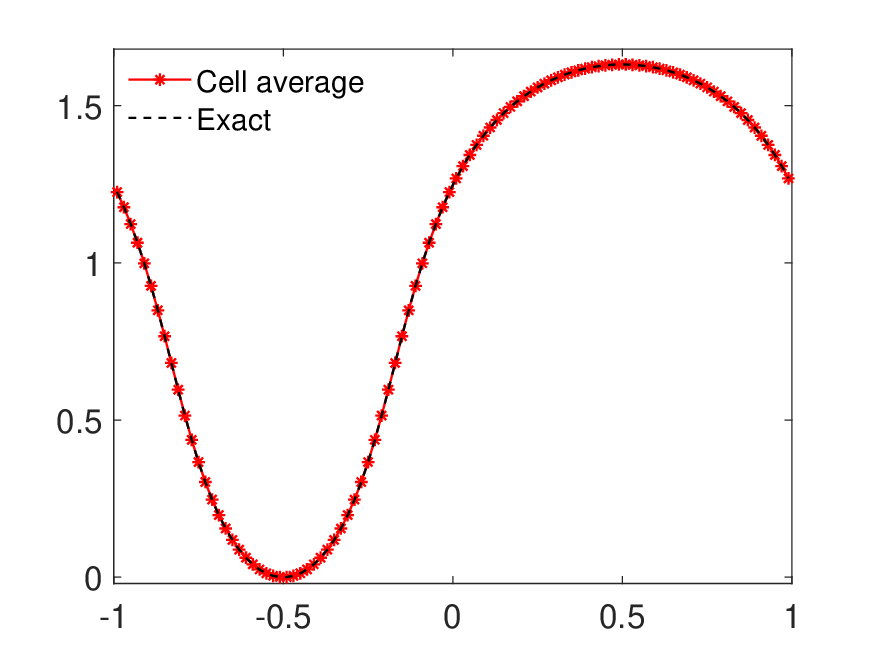}}\hspace*{0.15cm}
\subfigure[cell average of $v$]{\includegraphics[trim=0.8cm 0.4cm 0.9cm 0.5cm,clip,width=4.2cm]{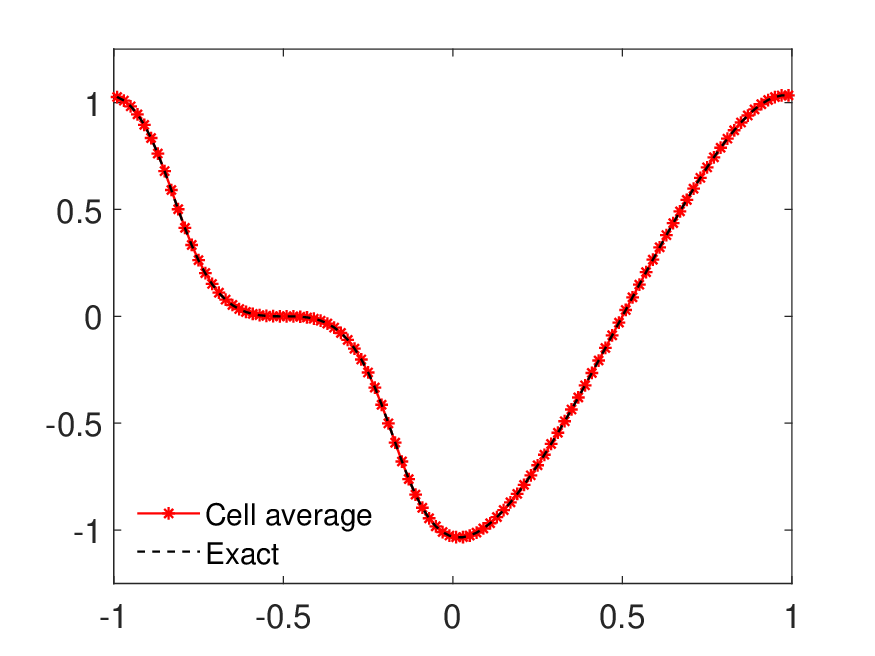}}\hspace*{0.15cm}
\subfigure[cell average of $p$]{\includegraphics[trim=0.8cm 0.4cm 0.9cm 0.5cm,clip,width=4.2cm]{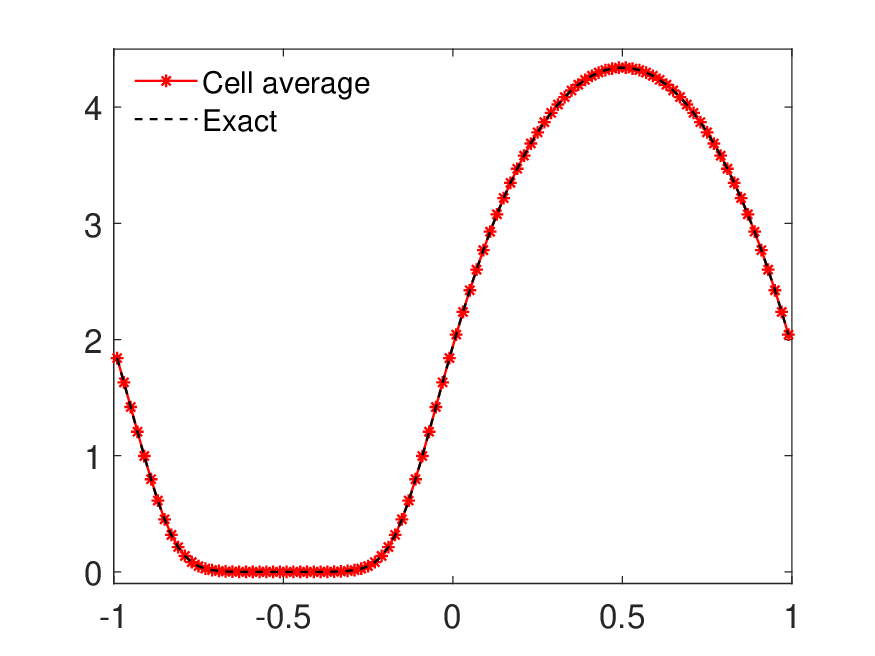}}}
\vskip5pt
\centerline{\subfigure[point value of $\rho$]{\includegraphics[trim=0.8cm 0.4cm 0.9cm 0.5cm,clip,width=4.2cm]{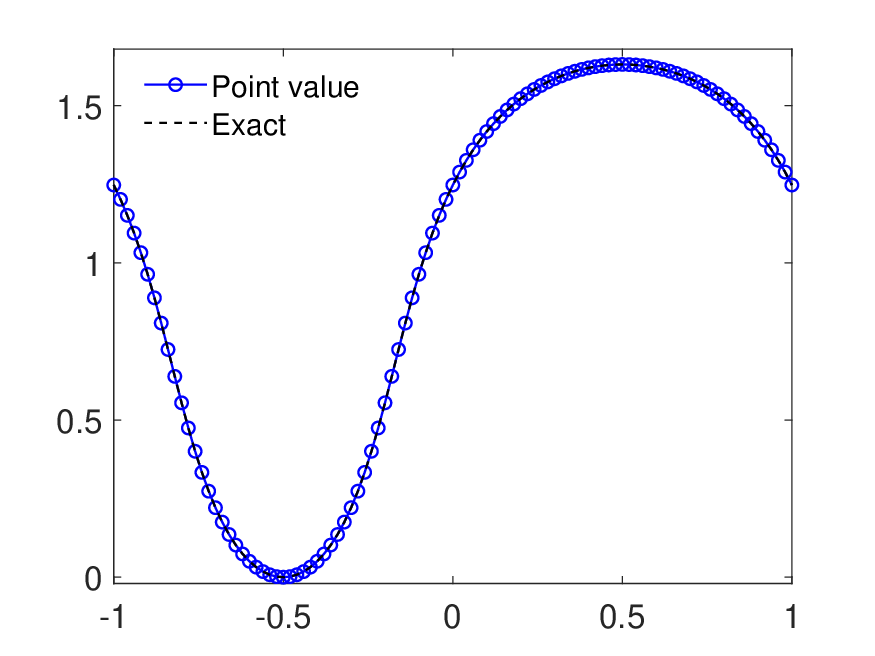}}\hspace*{0.15cm}
\subfigure[point values of $v$]{\includegraphics[trim=0.8cm 0.4cm 0.9cm 0.5cm,clip,width=4.2cm]{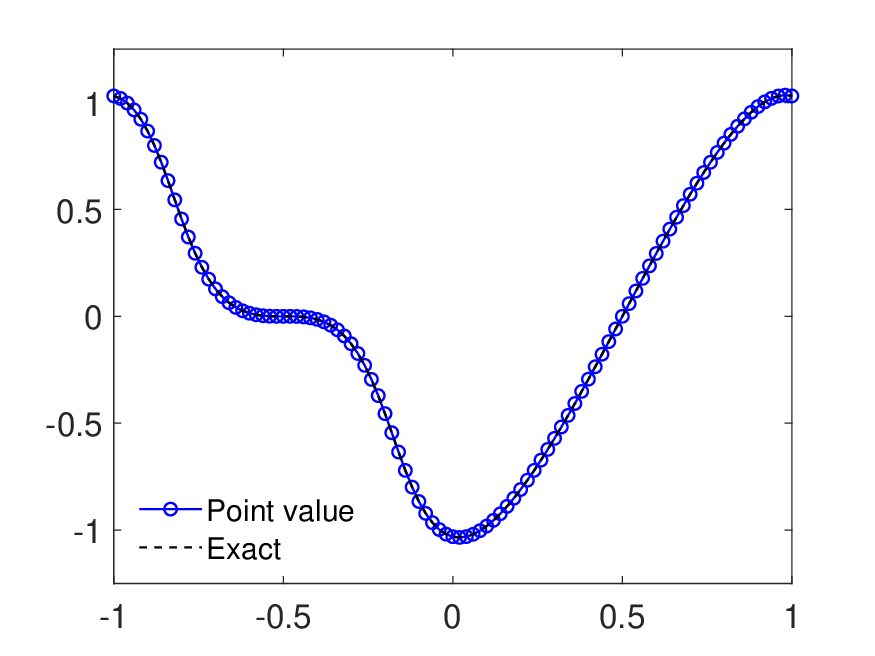}}\hspace*{0.15cm}
\subfigure[point value of $p$]{\includegraphics[trim=0.8cm 0.4cm 0.9cm 0.5cm,clip,width=4.2cm]{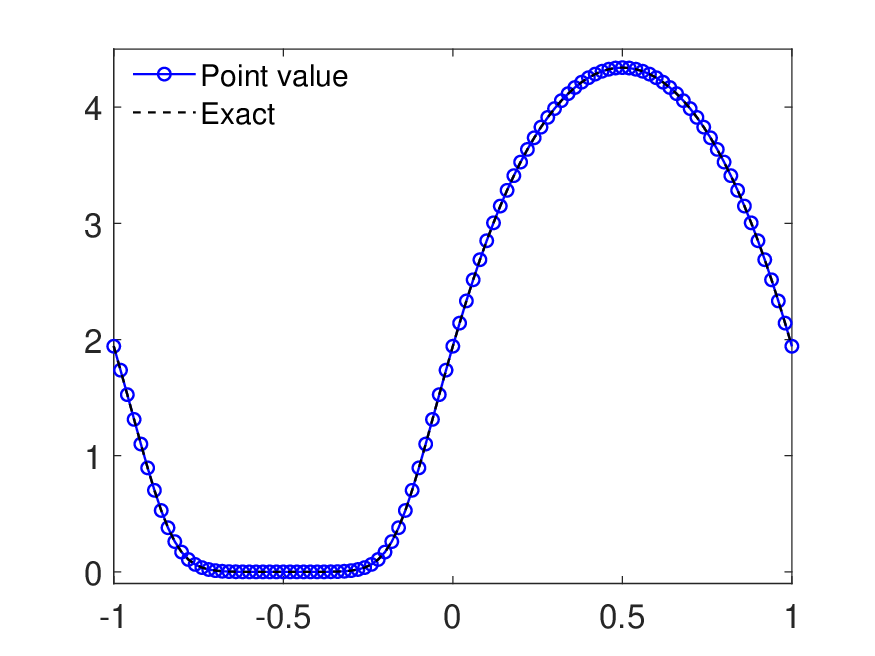}}}
\caption{\sf Example 4 (Isentropic): Solutions computed by BP PAMPA scheme {on a uniform mesh with 100 cells}. Top: cell average; Bottom: point value.\label{fig:Isen}}
\end{figure}

This example is challenging, as initially we have $\rho\left(-\frac{1}{2}\right) = 5.0 \times 10^{-6}$ and $p\left(-\frac{1}{2}\right) = 1.25 \times 10^{-16}$. With density and pressure values so close to zero, any numerical scheme lacking positivity preservation would fail. We test this example to examine the accuracy and robustness of the proposed BP PAMPA scheme, conducting the simulation until the final time $t = 0.1$. The solution, computed on a uniform mesh with 100 cells, is shown in Figure \ref{fig:Isen}, where we observe that density and pressure remain positive, even as they approach zero around $x = -\frac{1}{2}$.

Table \ref{tab1} gives the numerical errors and corresponding convergence rates. We also have looked for the degradation of accuracy using the BP procedure, as the solution has been locally corrected to preserve positivity. However, without BP limiting, the original third-order scheme quickly blows up: the pressure is so low in this example that it soon becomes negative without specific stabilization.
\begin{table}[ht!]
\caption{\label{tab1} Example 4: Errors and convergence rates in density for the isentropic problem.}
\begin{center}
\begin{tabular}{|c||c|c|c|c|c|c|}
\hline
$\dx$& $L^1$-error &rate&  $L^2$-error& rate&$L^\infty$-error& rate \\
\hline
$4\;10^{-2}$     &  $5.818\;10^{-4}$&- &  $6.581\;10^{-4}$& -     &        $1.390\;10{-3}$  & - \\
$2\;10^{-2}$ &  $8.056\;10^{-5}$&$2.85$&  $9.616\;10^{-5}$& $2.77$&       $2.212\;10^{-4}$  & $2.65$\\
$1\;10^{-2}$ &  $1.126\;10^{-5}$& $2.86$&  $1.363\;10^{-5}$& $2.82$&       $3.038\;10^{-5}$  & $2.84$\\
$5\;10^{-3}$ &  $1.602\;10^{-6}$& $2.82$&  $1.800\;10^{-6}$& $2.92$&       $3.836\;10^{-6}$  & $2.99$\\
$2.5\;10^{-3}$ &  $4.828\;10^{-7}$& $1.73$&  $3.780\;10^{-7}$& $2.25$&       $7.396\;10^{-7}$  & $2.37$\\
\hline
\end{tabular}
\end{center}
\end{table}

\subsubsection*{Example 5--- Sod problem}
\begin{figure}[ht!]
\centerline{\subfigure[cell average of $\rho$]{\includegraphics[trim=0.8cm 0.4cm 0.9cm 0.5cm,clip,width=4.2cm]{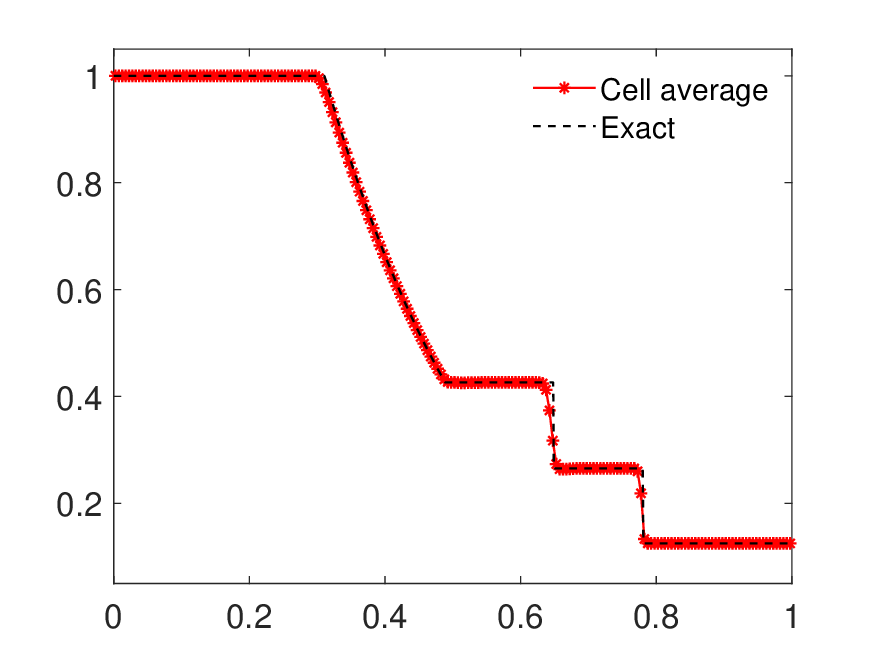}}\hspace*{0.15cm}
\subfigure[cell average of $v$]{\includegraphics[trim=0.8cm 0.4cm 0.9cm 0.5cm,clip,width=4.2cm]{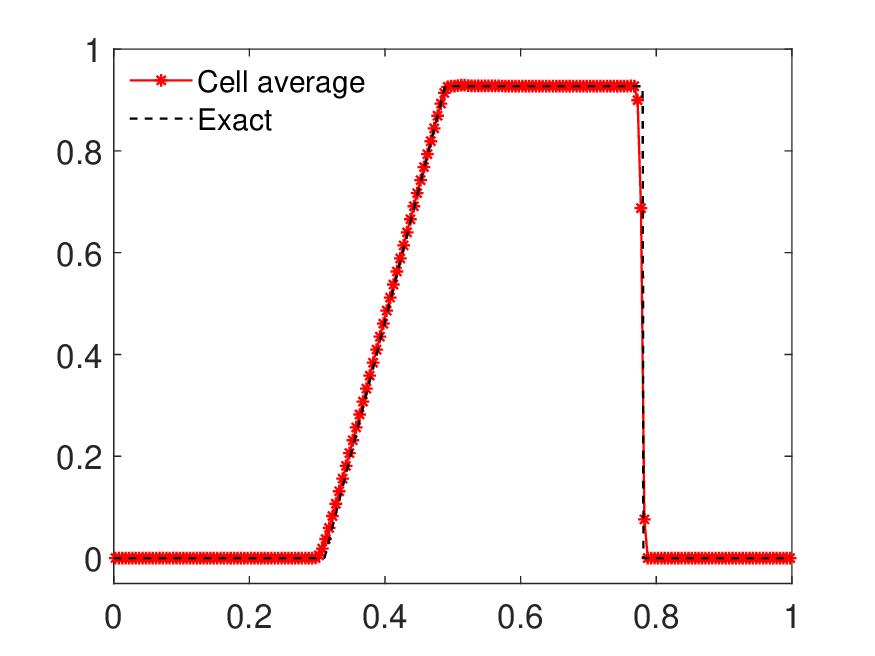}}\hspace*{0.15cm}
\subfigure[cell average of $p$]{\includegraphics[trim=0.8cm 0.4cm 0.9cm 0.5cm,clip,width=4.2cm]{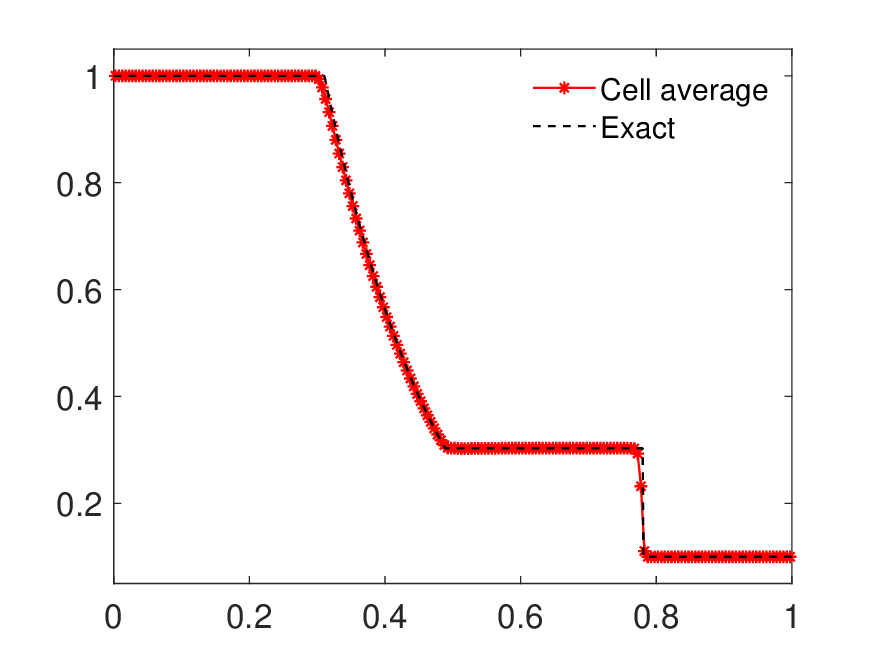}}}
\vskip5pt
\centerline{\subfigure[point value of $\rho$]{\includegraphics[trim=0.8cm 0.4cm 0.9cm 0.5cm,clip,width=4.2cm]{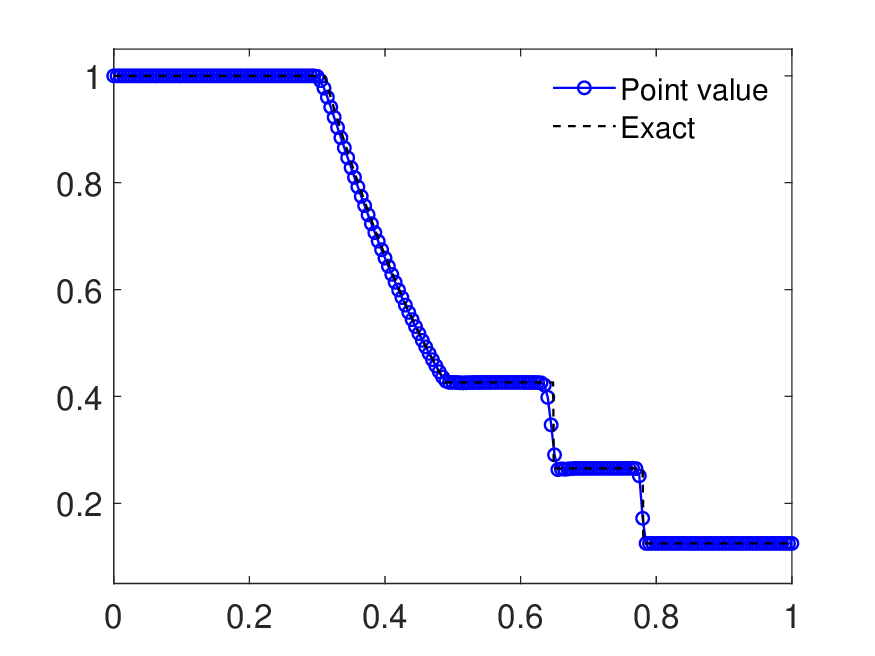}}\hspace*{0.15cm}
\subfigure[point values of $v$]{\includegraphics[trim=0.8cm 0.4cm 0.9cm 0.5cm,clip,width=4.2cm]{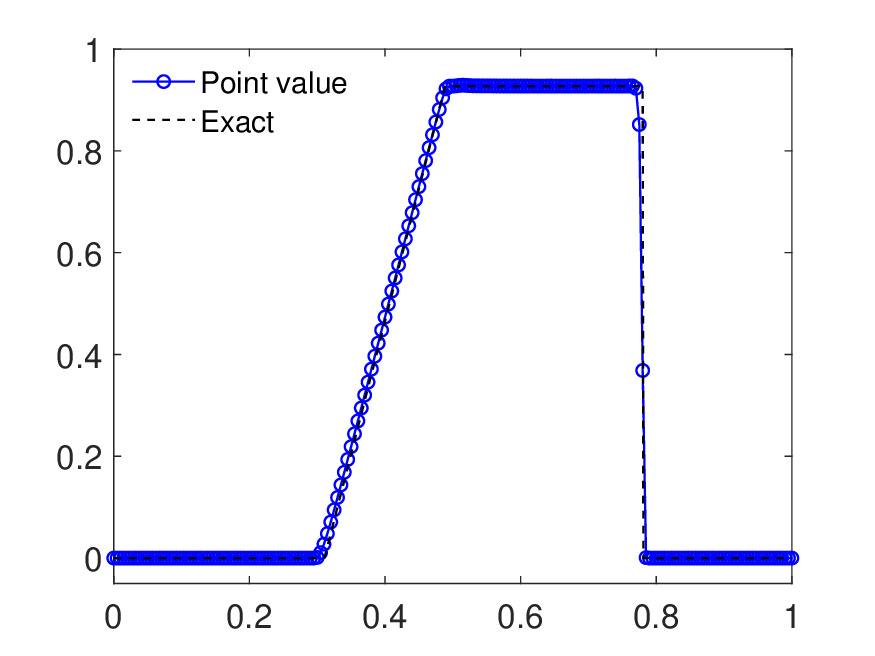}}\hspace*{0.15cm}
\subfigure[point value of $p$]{\includegraphics[trim=0.8cm 0.4cm 0.9cm 0.5cm,clip,width=4.2cm]{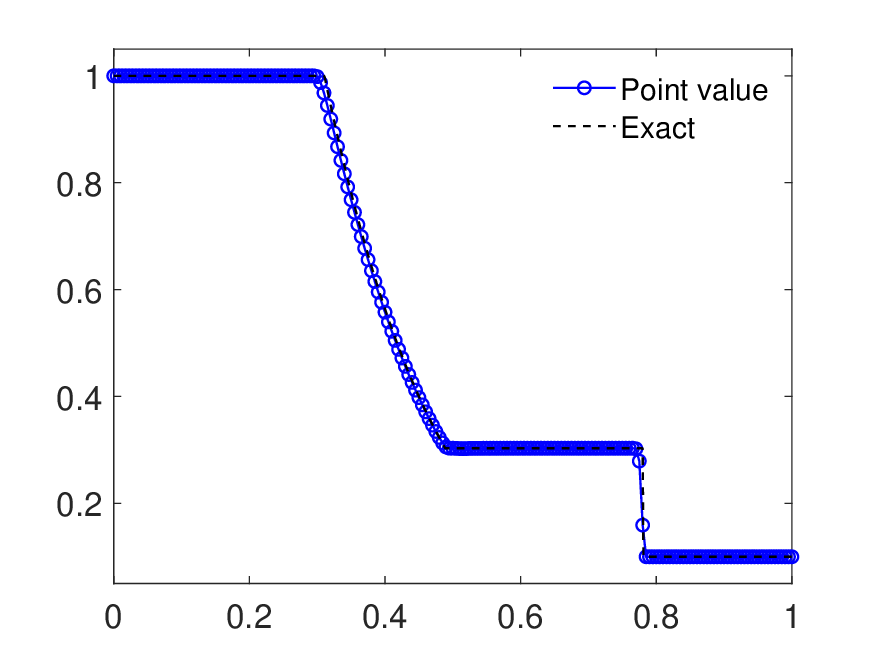}}}
\caption{\sf Example 5 (Sod): Solutions computed by BP PAMPA scheme {on a uniform mesh with 200 cells}. Top: cell average; Bottom: point value.\label{fig:Sod}}
\end{figure}

In the fifth example, we consider the Sod case defined in $[0,1]$ with the initial data given by
\begin{equation*}
    (\rho, v, p)=\left\{\begin{array}{ll}
         (1, 0, 1) &\mbox{if}~x<0.5,\\
         (0.125, 0, 0.1)&\mbox{else}.
    \end{array}\right.
\end{equation*}
The solution is computed until the final time $t=0.16$ on a uniform mesh with $200$ cells. The obtained results are displayed in Figure \ref{fig:Sod}. The proposed BP PAMPA scheme produces a correct solution with good resolution and performs well at the contact discontinuity.

\subsubsection*{Example 6---Shu-Osher problem}
\begin{figure}[ht!]
\centerline{\subfigure[cell average of $\rho$]{\includegraphics[trim=0.8cm 0.4cm 0.9cm 0.5cm,clip,width=4.2cm]{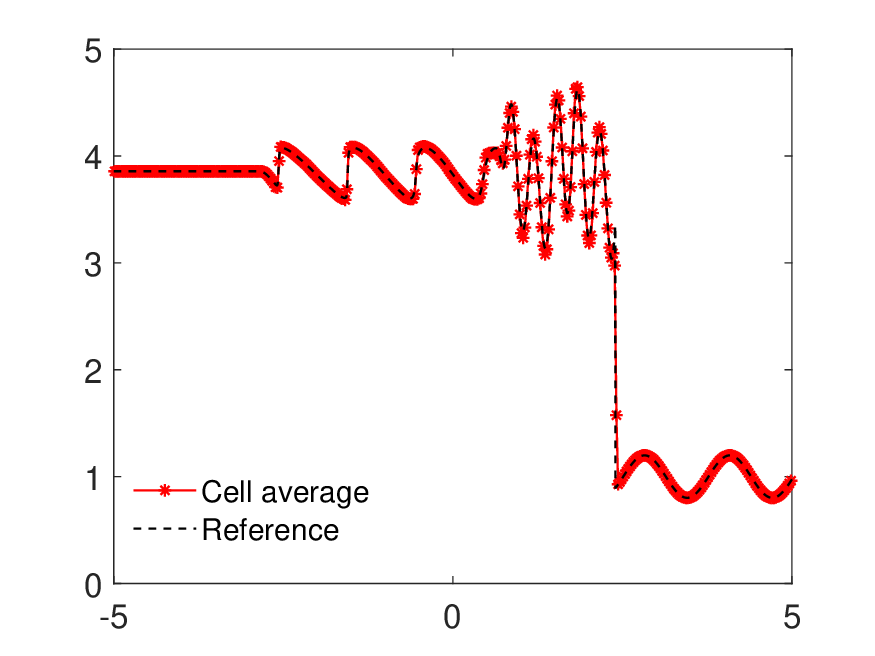}}\hspace*{0.15cm}
\subfigure[cell average of $v$]{\includegraphics[trim=0.8cm 0.4cm 0.9cm 0.5cm,clip,width=4.2cm]{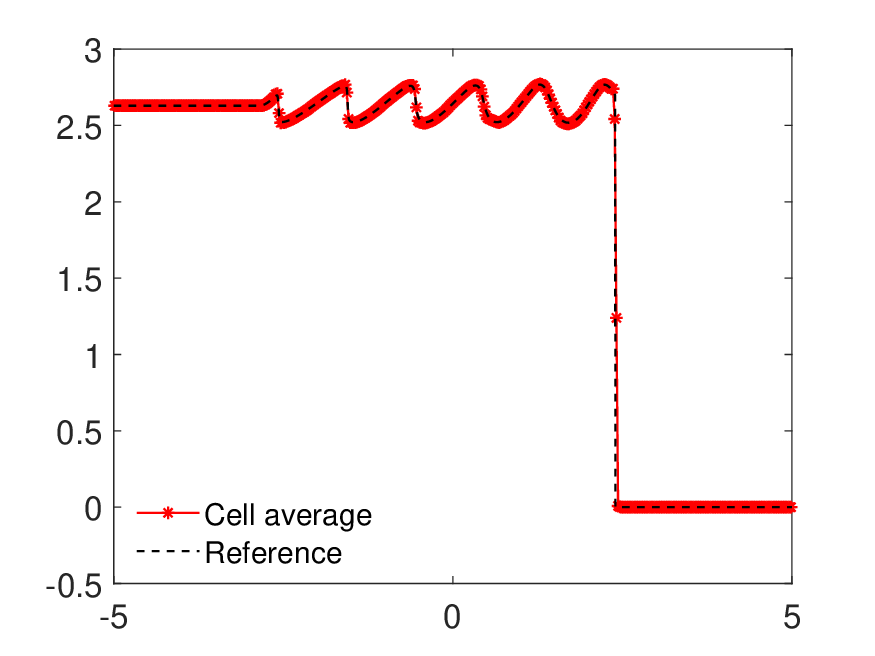}}\hspace*{0.15cm}
\subfigure[cell average of $p$]{\includegraphics[trim=0.8cm 0.4cm 0.9cm 0.5cm,clip,width=4.2cm]{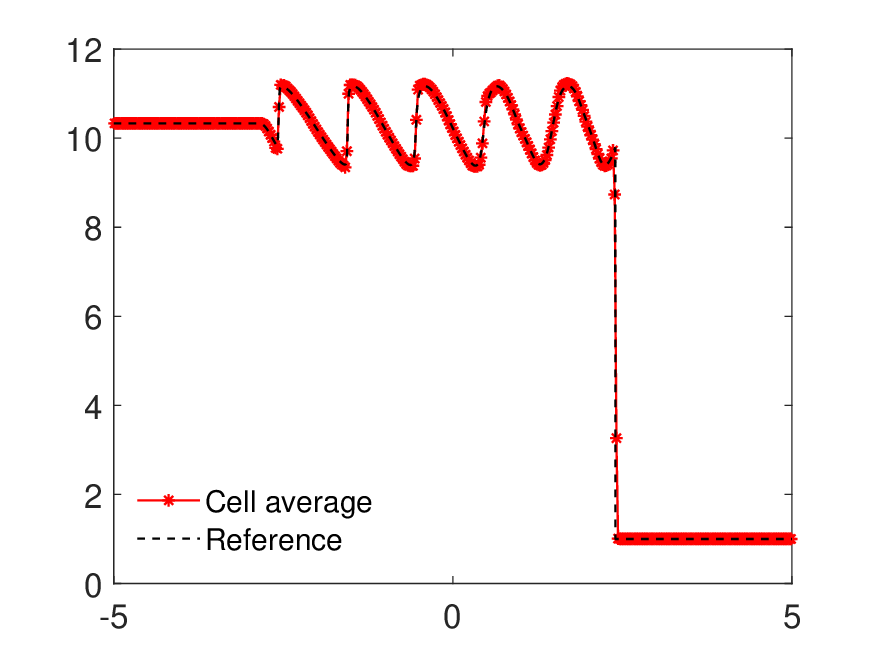}}}
\vskip5pt
\centerline{\subfigure[point value of $\rho$]{\includegraphics[trim=0.8cm 0.4cm 0.9cm 0.5cm,clip,width=4.2cm]{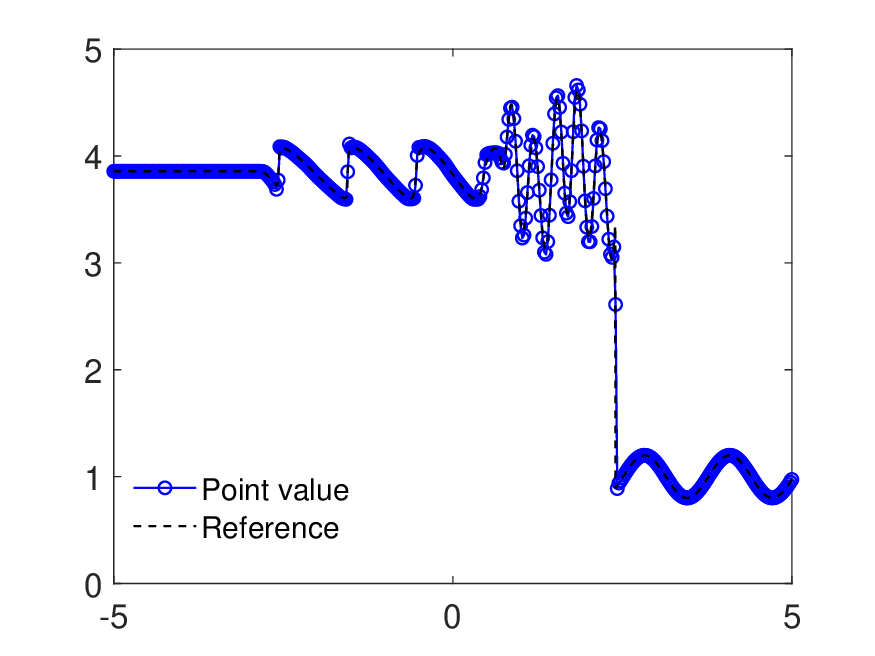}}\hspace*{0.15cm}
\subfigure[point values of $v$]{\includegraphics[trim=0.8cm 0.4cm 0.9cm 0.5cm,clip,width=4.2cm]{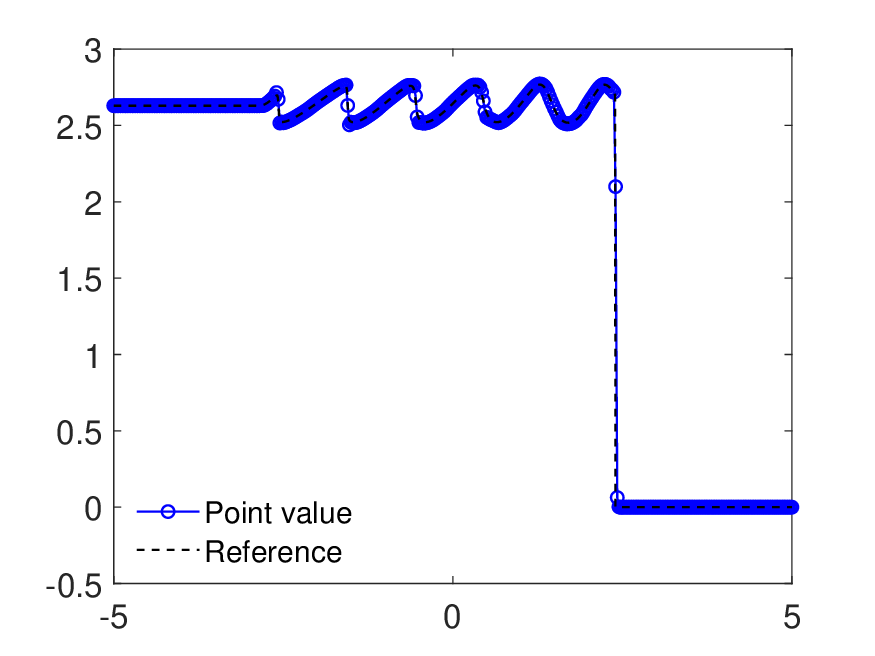}}\hspace*{0.15cm}
\subfigure[point value of $p$]{\includegraphics[trim=0.8cm 0.4cm 0.9cm 0.5cm,clip,width=4.2cm]{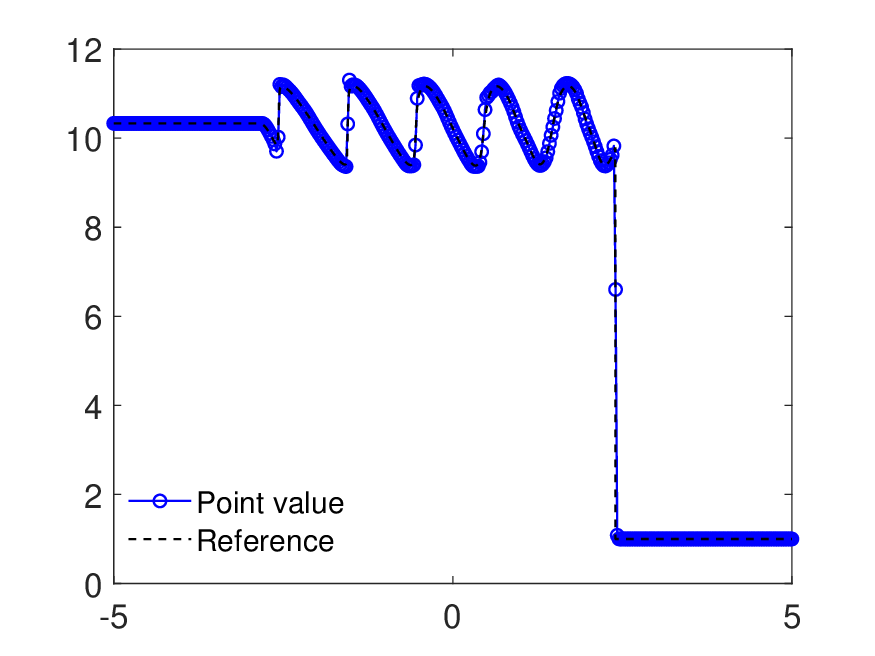}}}
\caption{\sf Example 6 (Shu-Osher): Solutions computed by BP PAMPA scheme {on a uniform mesh with 400 cells. The reference solution is computed by the first-order local Lax--Friedrichs scheme using 40,000 cells}. Top: cell average; Bottom: point value.\label{fig:Shu-Osher}}
\end{figure}

In the sixth example, we consider the shock and sine wave interaction problem from \cite{Shu1989}. The initial data is
\begin{equation*}
    (\rho, v, p)=\left\{\begin{array}{ll}
         (3.857143, 2.629369, 10.33333333333) &\mbox{if}~x<-4,\\
         (1+0.2\sin(5x), 0, 1)&\mbox{else},
    \end{array}\right.
\end{equation*}
prescribed in the computational domain $[-5,5]$. We compute the solution using the BP PAMPA scheme until the final time $t=1.8$ on a uniform mesh with 400 cells. The results are shown in Figure \ref{fig:Shu-Osher}, along with a reference solution computed by the first-order local Lax--Friedrichs scheme using 40,000 uniform cells. We observe that the BP PAMPA scheme accurately captures the interaction of sine waves and the right-moving shock, consistently providing good results without spurious oscillations.

\subsubsection*{Example 7---LeBlanc problem}
\begin{figure}[ht!]
\centerline{\subfigure[cell average of $\rho$]{\includegraphics[trim=0.8cm 0.4cm 0.9cm 0.5cm,clip,width=4.2cm]{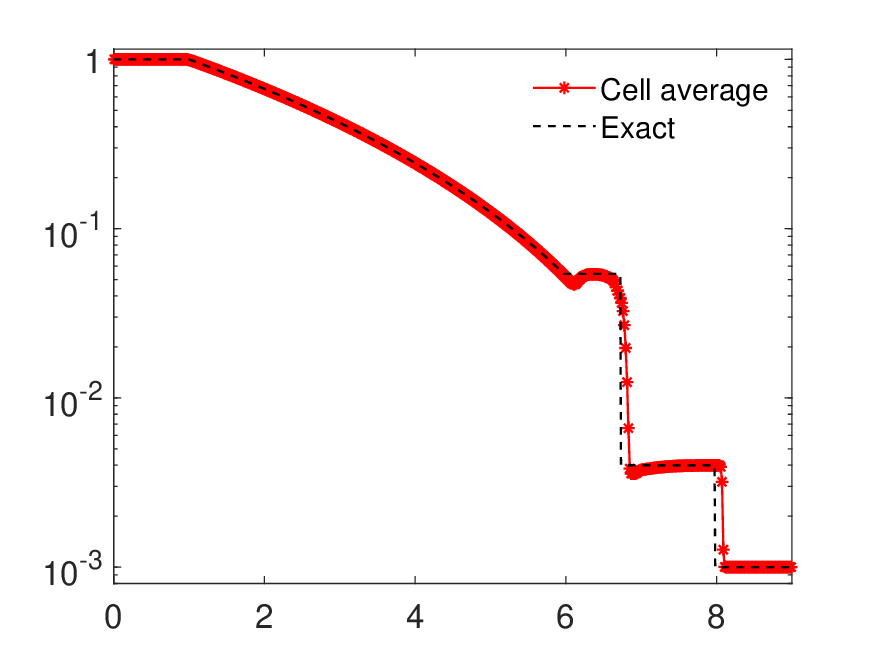}}\hspace*{0.15cm}
\subfigure[cell average of $v$]{\includegraphics[trim=0.8cm 0.4cm 0.9cm 0.5cm,clip,width=4.2cm]{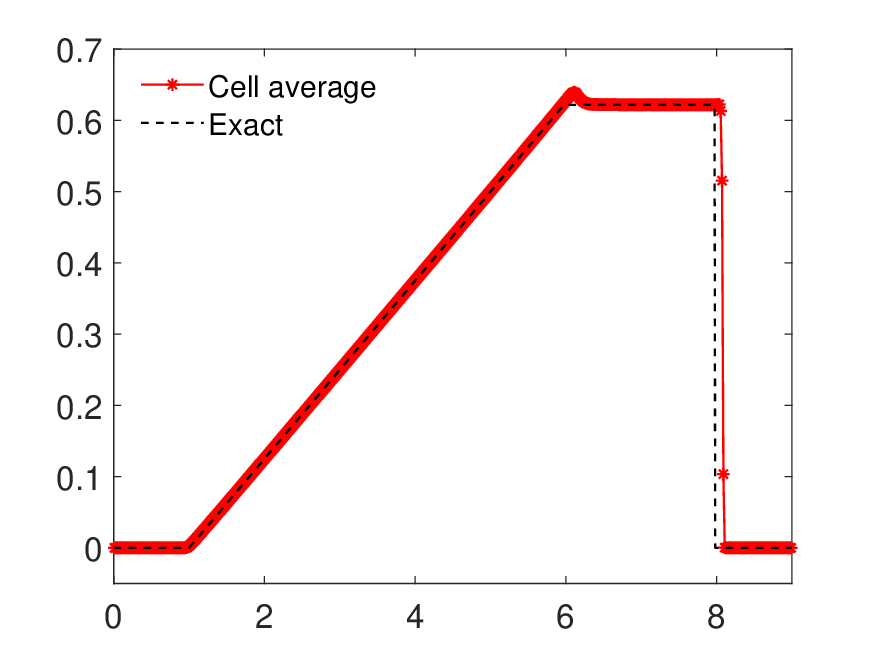}}\hspace*{0.15cm}
\subfigure[cell average of $p$]{\includegraphics[trim=0.8cm 0.4cm 0.9cm 0.5cm,clip,width=4.2cm]{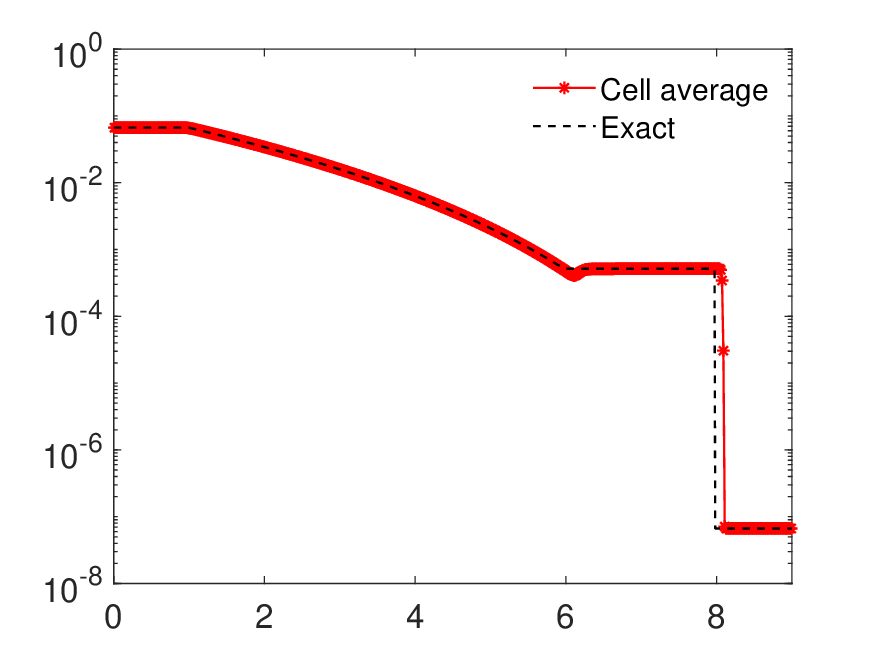}}}
\vskip5pt
\centerline{\subfigure[point value of $\rho$]{\includegraphics[trim=0.8cm 0.4cm 0.9cm 0.5cm,clip,width=4.2cm]{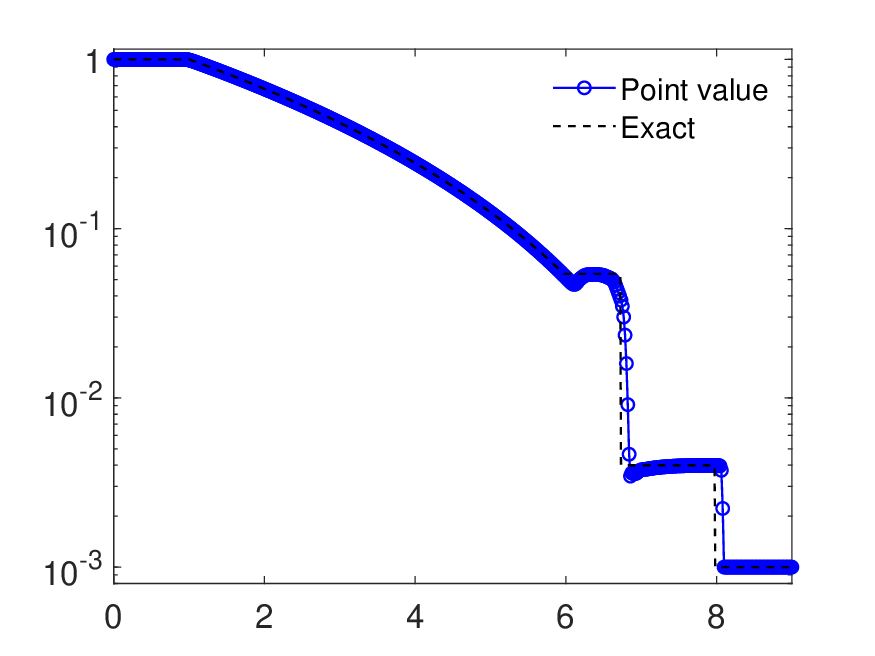}}\hspace*{0.15cm}
\subfigure[point values of $v$]{\includegraphics[trim=0.8cm 0.4cm 0.9cm 0.5cm,clip,width=4.2cm]{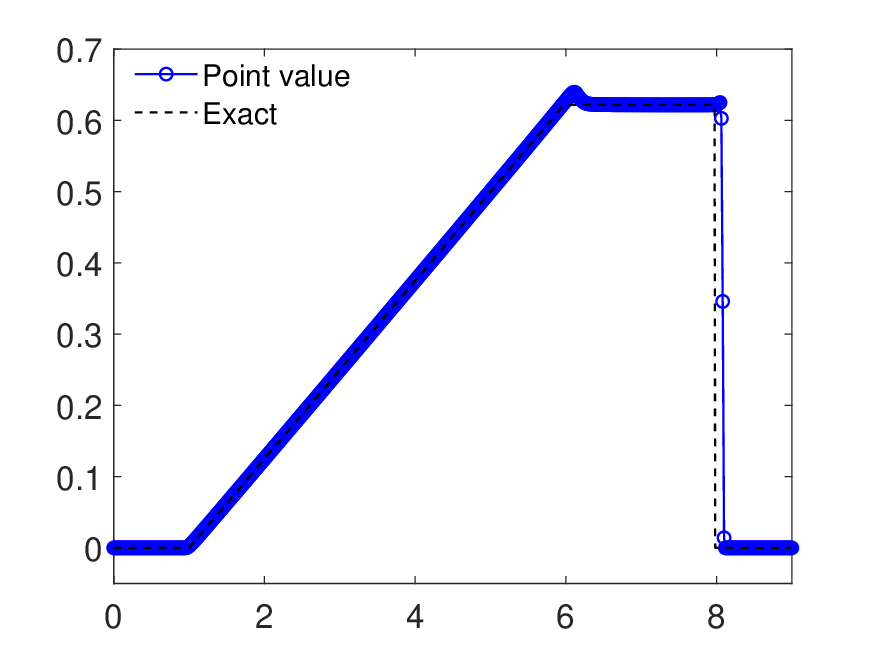}}\hspace*{0.15cm}
\subfigure[point value of $p$]{\includegraphics[trim=0.8cm 0.4cm 0.9cm 0.5cm,clip,width=4.2cm]{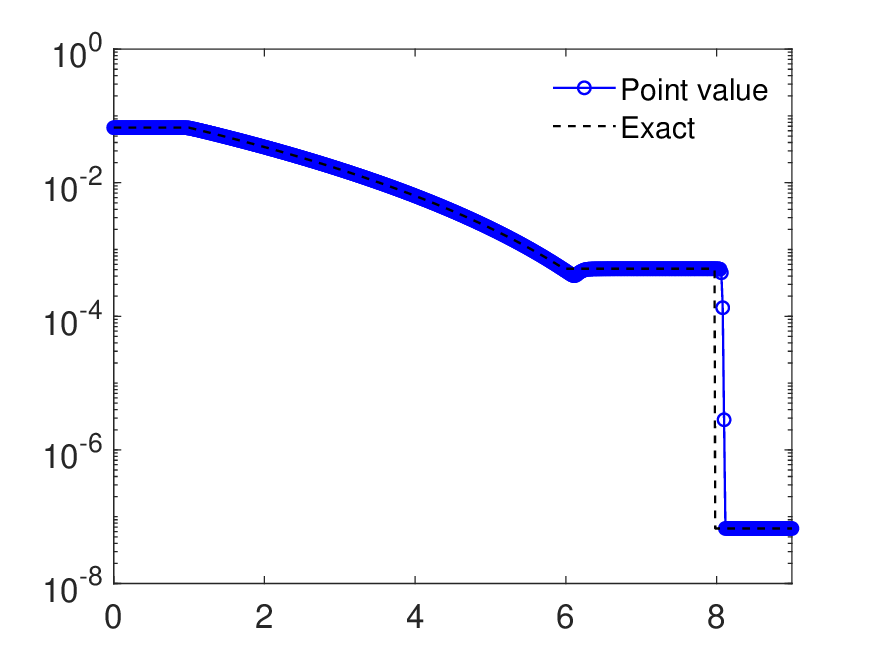}}}
\caption{\sf Example 7 (LeBlanc): Solutions computed by BP PAMPA scheme {on a uniform mesh with 500 cells}. Top: cell average; Bottom: point value.\label{fig:LeBlanc}}
\end{figure}

In the seventh example, we study the LeBlanc problem on the computational domain $[0,9]$. The initial data is defined as
\begin{equation*}
    (\rho, v, p)=\left\{\begin{array}{ll}
         (1, 0, 0.1(\gamma-1)) &\mbox{if}~x \leq 3,\\
         (0.001, 0, 10^{-7}(\gamma-1)) &\mbox{else}.
    \end{array}\right.
\end{equation*}
This is a very strong shock tube problem, and without BP limiters, it is challenging or even impossible to run simulations above first-order due to the risk of negative density or pressure. Figure \ref{fig:LeBlanc} shows the solutions computed by the BP PAMPA scheme at the final time $t=6$, using a uniform mesh with 500 cells. As observed, a strong shock wave propagates from the high-pressure region on the left to the low-density region on the right, while a rarefaction wave moves to the left.

\begin{figure}[ht!]
\centerline{\subfigure[cell average of $\rho$]{\includegraphics[trim=0.8cm 0.4cm 0.9cm 0.5cm,clip,width=4.2cm]{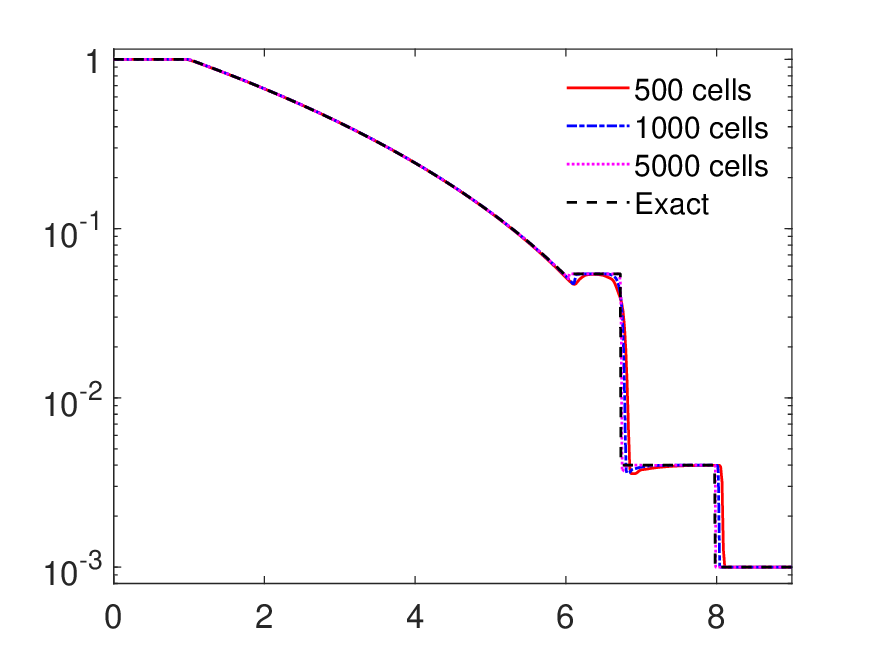}}\hspace*{0.15cm}
\subfigure[zoom for $6\leq x\leq9$]{\includegraphics[trim=0.8cm 0.4cm 0.9cm 0.5cm,clip,width=4.2cm]{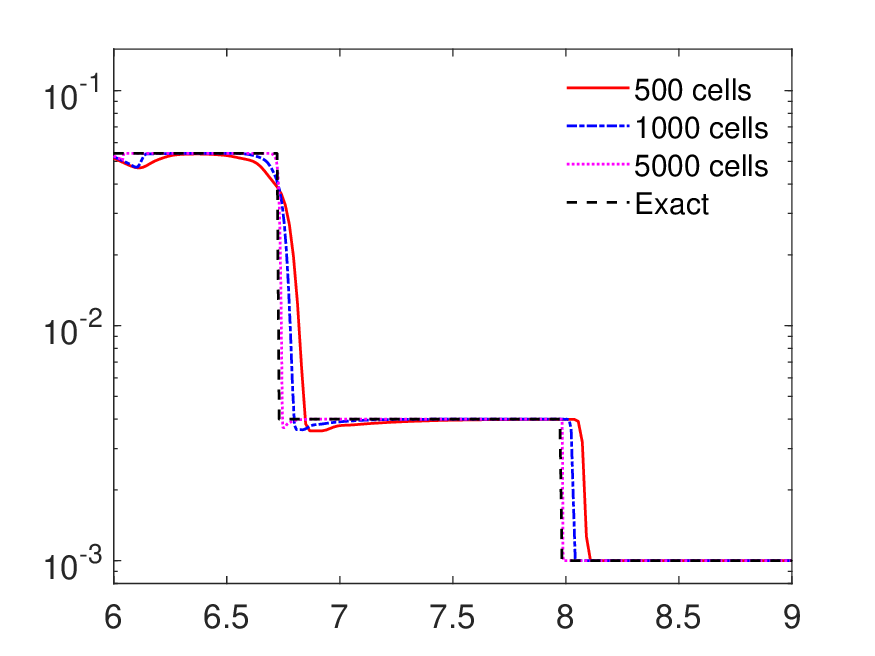}}\hspace*{0.15cm}
\subfigure[zoom for $7.5\leq x\leq8.5$ ]{\includegraphics[trim=0.8cm 0.4cm 0.9cm 0.5cm,clip,width=4.2cm]{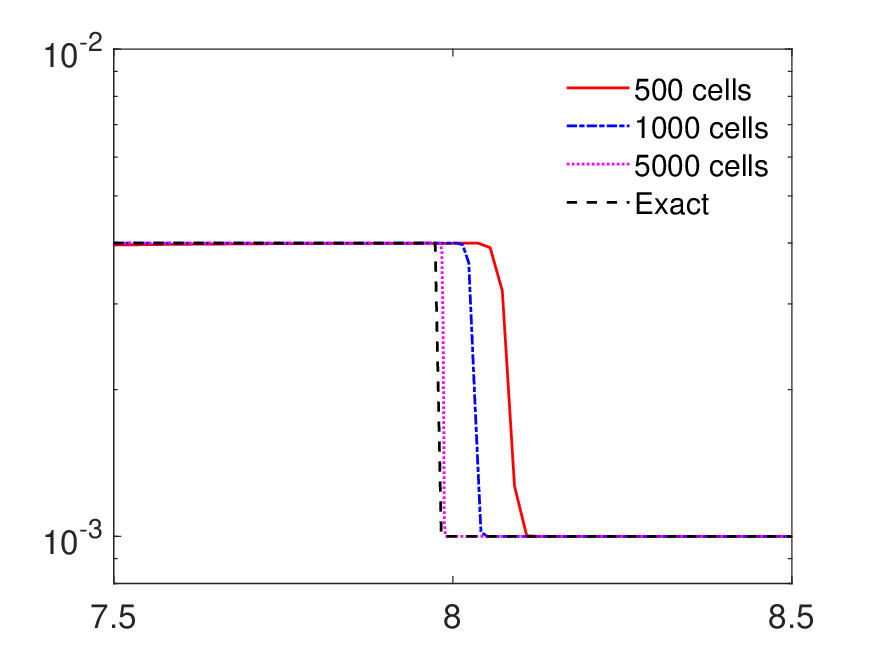}}}
\caption{\sf Example 7 (LeBlanc): Convergence study on the average value of density.\label{fig:LeBlanc2}}
\end{figure}

It is known that at time $t=6$, the shock wave should be located at $x=8$. Given the extreme conditions, accurately capturing the shock wave position is generally challenging. We therefore perform a convergence study, with results shown in Figure \ref{fig:LeBlanc2}. As seen, the numerical solution converges to the exact one, and the results are as expected.

\subsubsection*{Example 8---Double rarefaction problem}
\begin{figure}[ht!]
\centerline{\subfigure[cell average of $\rho$]{\includegraphics[trim=0.8cm 0.4cm 0.9cm 0.5cm,clip,width=4.2cm]{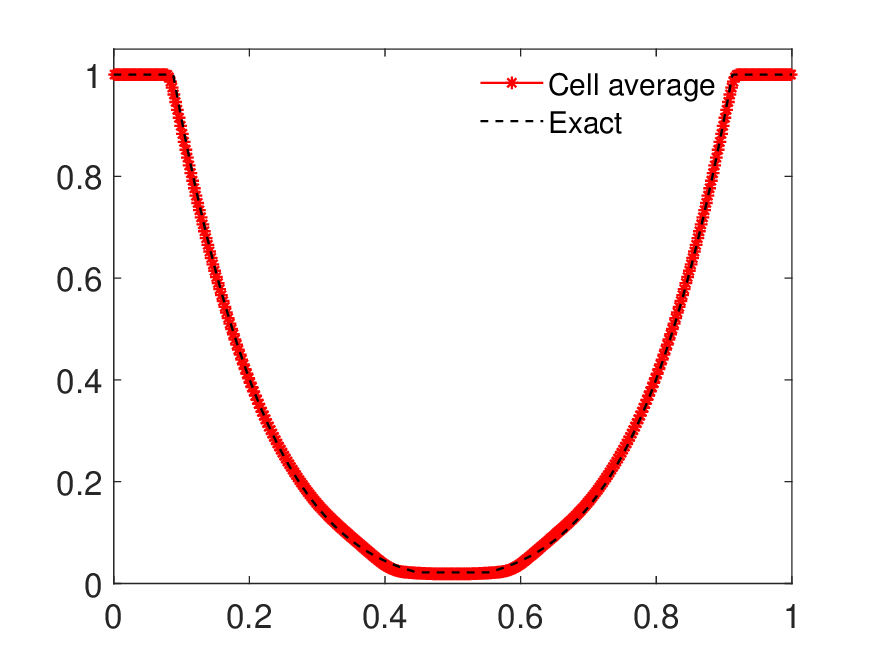}}\hspace*{0.15cm}
\subfigure[cell average of $v$]{\includegraphics[trim=0.8cm 0.4cm 0.9cm 0.5cm,clip,width=4.2cm]{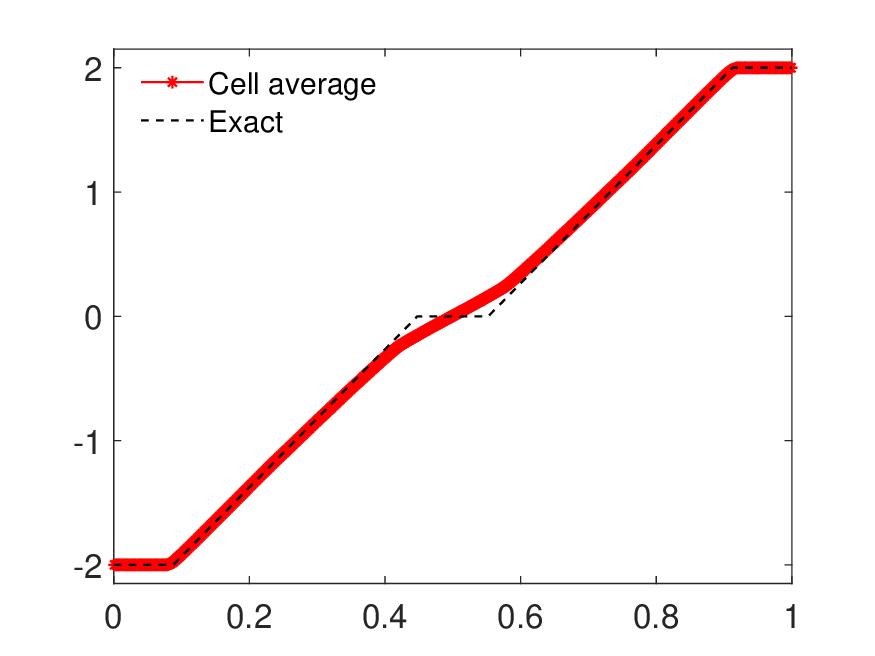}}\hspace*{0.15cm}
\subfigure[cell average of $p$]{\includegraphics[trim=0.8cm 0.4cm 0.9cm 0.5cm,clip,width=4.2cm]{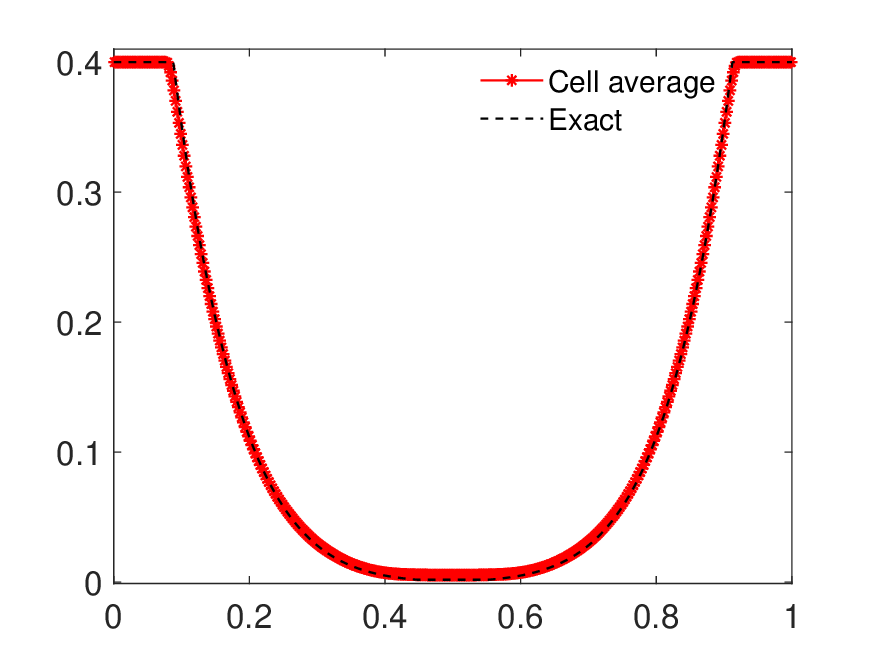}}}
\vskip5pt
\centerline{\subfigure[point value of $\rho$]{\includegraphics[trim=0.8cm 0.4cm 0.9cm 0.5cm,clip,width=4.2cm]{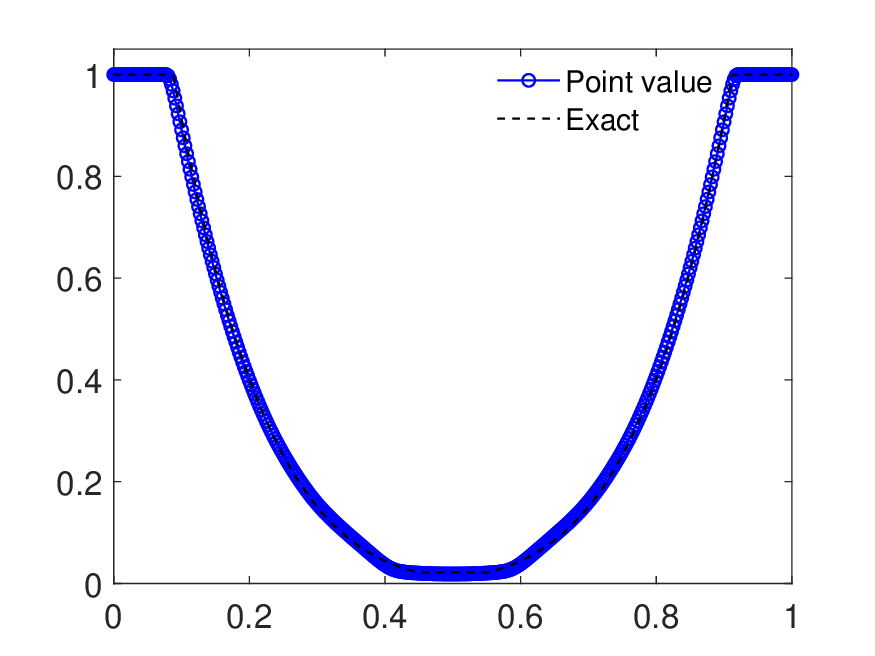}}\hspace*{0.15cm}
\subfigure[point values of $v$]{\includegraphics[trim=0.8cm 0.4cm 0.9cm 0.5cm,clip,width=4.2cm]{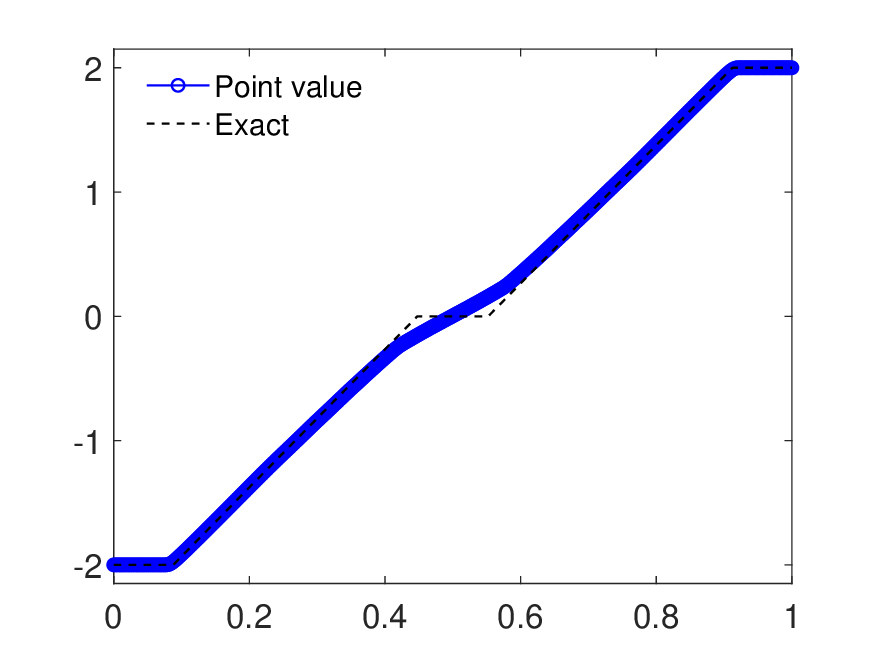}}\hspace*{0.15cm}
\subfigure[point value of $p$]{\includegraphics[trim=0.8cm 0.4cm 0.9cm 0.5cm,clip,width=4.2cm]{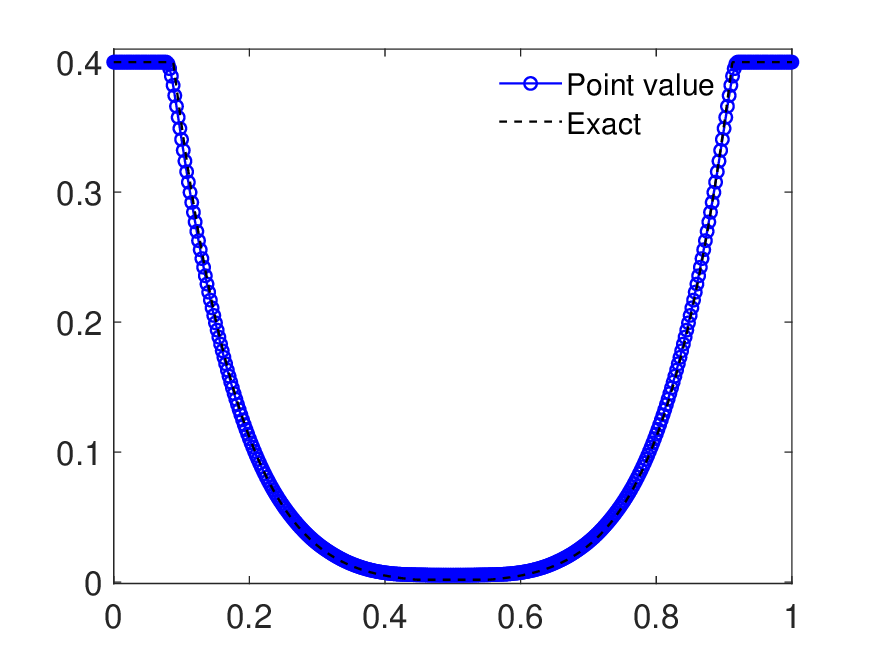}}}
\caption{\sf Example 8 (double rarefaction): Solutions computed by BP PAMPA scheme {on a uniform mesh with 400 cells}. Top: cell average; Bottom: point value.\label{fig:double-rare}}
\end{figure}

In the eighth example, we study the so-called ``1-2-3'' problem with the following initial data
\begin{equation*}
    (\rho, v, p)=\left\{\begin{array}{ll}
         (1, -2, 0.4)&\mbox{if}~x\leq0.5,\\
         (1, 2, 0.4)&\mbox{else},
    \end{array}\right.
\end{equation*}
prescribed in the computational domain $[0,1]$. The exact solution of this Riemann problem consists of two rarefaction waves with a near-vacuum region between them. We compute the solution until the final time $t=0.15$ on a uniform mesh with 400 cells and present the results in Figure \ref{fig:double-rare}. The proposed BP PAMPA scheme produces good, non-oscillatory results, showing two rarefaction waves and a near-vacuum region in the middle. 

\subsubsection*{Example 9---Sedov problem}
In the final example, we consider the Sedov problem \cite{Sedov,vonNeuman}, which features very low density with strong shocks. Initialization is done as follows (see \cite{MaireHDR}):
\begin{itemize}
    \item The mesh points are $x_i=i\Delta x$, with $\Delta x=\frac{2}{N}$, $N=401$ and $-N\leq i\leq N$;
    \item For all $i$, the point value data are $\big(\rho_i, u_i, p_i\big)=\big ( 1,0,\frac{e_{\min}}{(\gamma-1)}\big )$ and $e_{\min}=10^{-12}$;
    \item For all cells $[x_i,x_{i+1}]$, $i\neq 0$, then $\big(\rho_i, \rho_i v_i, E_i\big)=\big(1,0,\rho_i e_{\min}\big)$;
    \item For the cell $[x_{0},x_1]=[-\frac{\Delta x}{2},\frac{\Delta x}{2}]$, we set $\big(\rho_i, \rho_i v_i, \rho_i \frac{e_{\max}}{\gamma-1}\big)=\big(1,0, \frac{\rho_i e_{\max}}{\gamma-1}\big)$ with $e_{\max}=\frac{0.538548}{\Delta x}$.
\end{itemize}
These conditions are set so that the strength of the shock is almost infinite,  and the initial condition corresponds to a Dirac mass or energy at $x=0$ in the limit of mesh refinement. The exact solution is computed following \cite{Kamm2000}.

\begin{figure}[ht!]
\centerline{\subfigure[cell average of $\rho$]{\includegraphics[trim=0.8cm 0.4cm 0.9cm 0.5cm,clip,width=4.2cm]{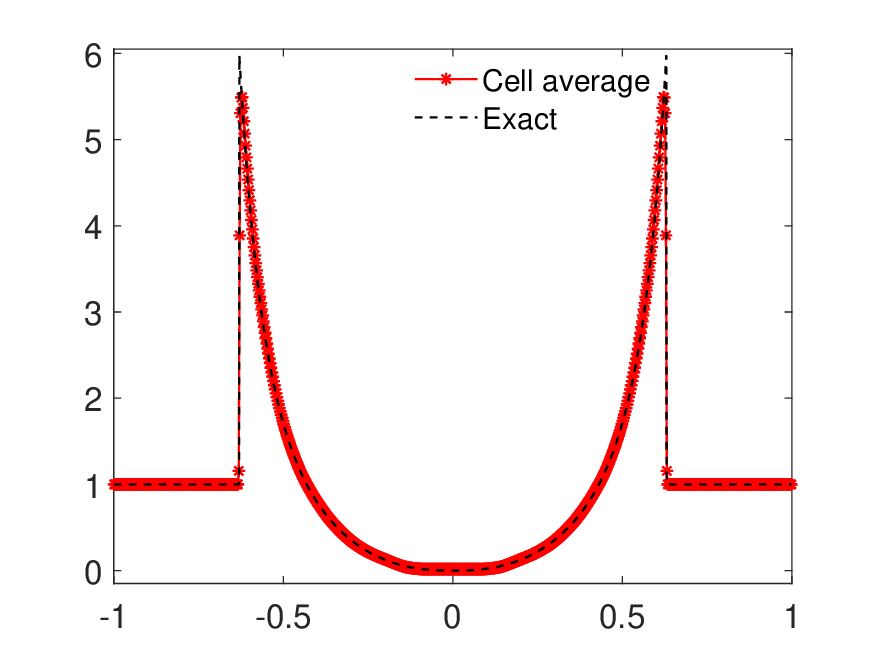}}\hspace*{0.15cm}
\subfigure[cell average of $v$]{\includegraphics[trim=0.8cm 0.4cm 0.9cm 0.5cm,clip,width=4.2cm]{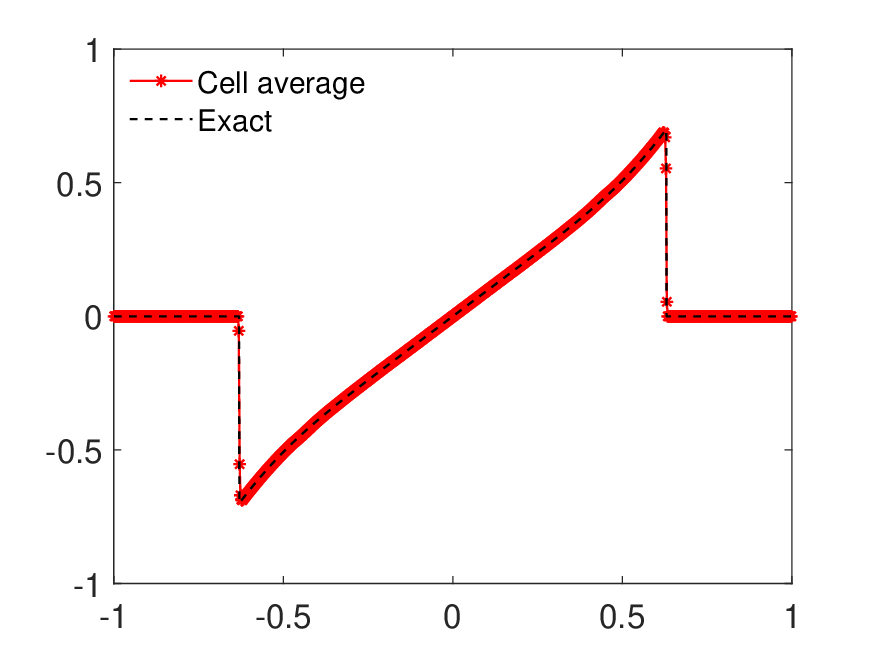}}\hspace*{0.15cm}
\subfigure[cell average of $p$]{\includegraphics[trim=0.8cm 0.4cm 0.9cm 0.5cm,clip,width=4.2cm]{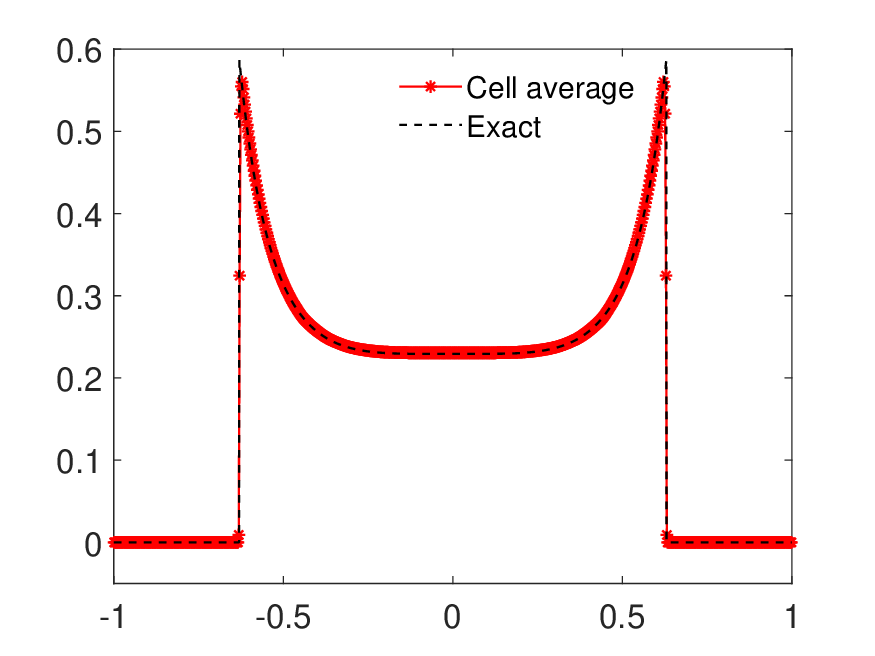}}}
\vskip5pt
\centerline{\subfigure[point value of $\rho$]{\includegraphics[trim=0.8cm 0.4cm 0.9cm 0.5cm,clip,width=4.2cm]{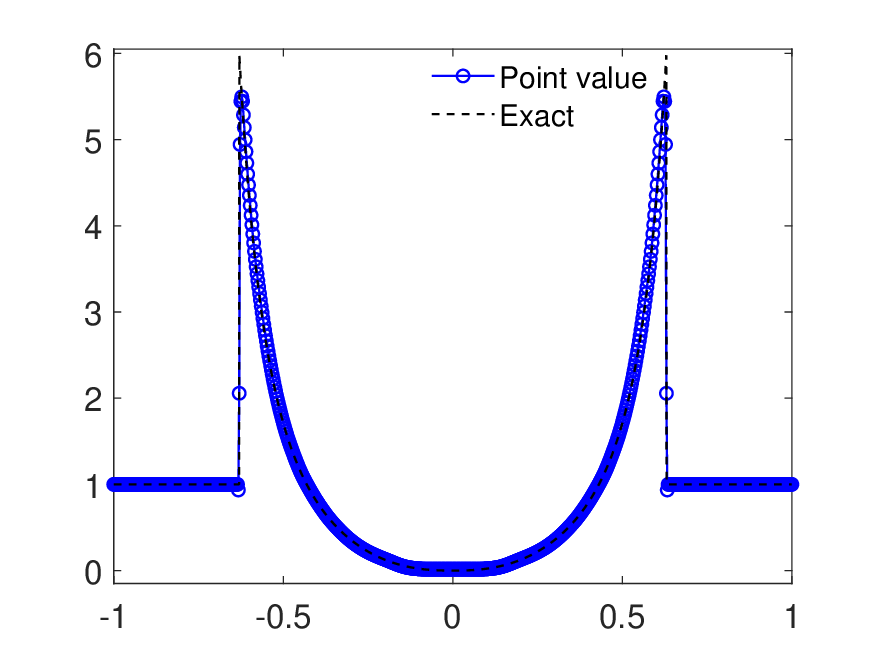}}\hspace*{0.15cm}
\subfigure[point values of $v$]{\includegraphics[trim=0.8cm 0.4cm 0.9cm 0.5cm,clip,width=4.2cm]{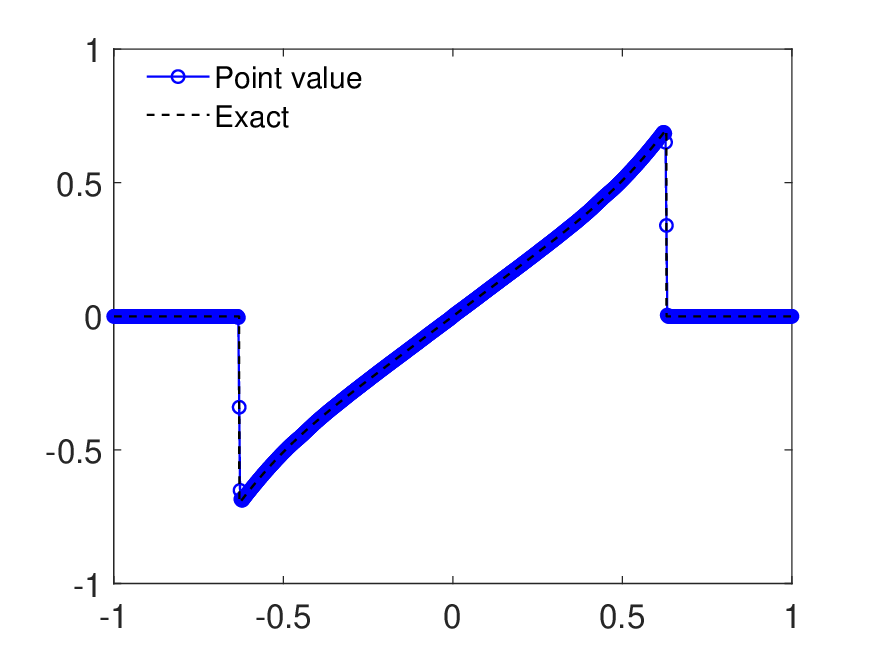}}\hspace*{0.15cm}
\subfigure[point value of $p$]{\includegraphics[trim=0.8cm 0.4cm 0.9cm 0.5cm,clip,width=4.2cm]{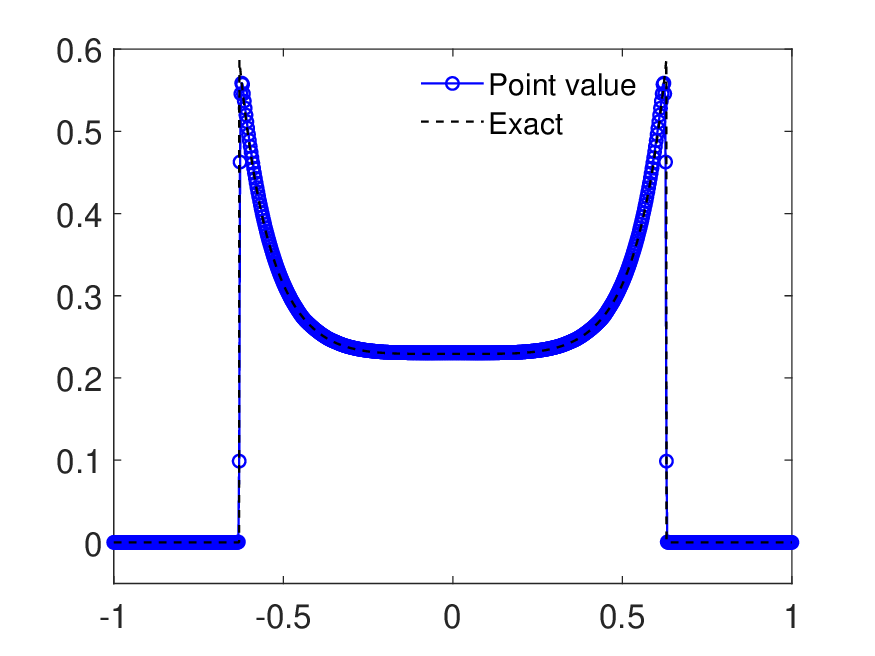}}}
\caption{\sf Example 9 (Sedov): Solutions computed by BP PAMPA scheme {on a uniform mesh with 401 cells}. Top: cell average; Bottom: point value.\label{fig:Sedov}}
\end{figure}

We compute the solution until the final time $t = 0.5$ and present the results in Figure \ref{fig:Sedov}. The BP PAMPA scheme produces satisfactory results, without any numerical defects.

\section{Conclusions}\label{sec5}
We have proposed a {BP} framework adapted to the PAMPA scheme, where the solution is represented by cell averages and point values. In this approach, the degrees of freedom are updated using a scheme that blends a first-order BP update with a high-order update. We have developed optimal blending parameters for both scalar and system cases. For hyperbolic systems, we employ the GQL approach of Wu and Shu, allowing us to propose explicit and optimal coefficients. This approach effectively preserves bounds and controls spurious oscillations for both the point values and cell averages, without requiring explicit solution reconstruction.

We intend to extend this work in several directions, such as higher-order approaches and multidimensional polygonal cells.

\appendix
\section{Minimization of $\varphi(\nu)$ knowing that $Q(\nu)> 0$}\label{appendix:A}
The aim of this appendix is to solve the following minimization problem
\begin{equation*}
   \min\limits_{\nu\in\mathbb{R}}\frac{Q(\nu)}{\vert P(\nu)\vert}, \quad Q(\nu)=a_1\nu^2-2b_1\nu+2c_1,\quad P(\nu)=a_2\nu^2-2b_2\nu+2c_2.
\end{equation*}

We assume that $\nu=\frac{p}{q}$ with $p,q\neq0$ and obtain
\begin{equation*}
    \frac{Q\big( p/q\big )}{P\big(p/q\big )}=\frac{a_1p^2-2b_1 pq+2c_1 q^2}{\vert a_2p^2-2b_2 pq+2c_2 q^2\vert}=\frac{\langle \mathbf u, B\mathbf u\rangle}{\vert \langle\mathbf u, A\mathbf u\rangle\vert},
\end{equation*}
where
\begin{equation*}
    \mathbf u=\begin{pmatrix} p\\q\end{pmatrix}, \quad A=\begin{pmatrix}a_2 & -b_2\\-b_2&2c_2\end{pmatrix}, \quad B=\begin{pmatrix}a_1 & -b_1\\-b_1&2c_1\end{pmatrix}.
\end{equation*}

Since $B$ is positive definite corresponding to positive eigenvalues, the minimization problem is equivalent to 
\begin{equation*}
 \min\limits_{\substack{\nu\in\mathbb{R}}}\frac{Q(\nu)}{\vert P(\nu)\vert}=\min\limits_{\substack{\mathbf u\in\mathbb{R}^2 \\\mathbf u\neq0}}\frac{\langle \mathbf u, B\mathbf u\rangle}{\vert\langle\mathbf u, A\mathbf u\rangle\vert}=\left(\max\limits_{\substack{\mathbf u\in\mathbb{R}^2\\\mathbf u\neq0}}{\frac{\vert\langle\mathbf u, A\mathbf u\rangle\vert}{\langle \mathbf u, B\mathbf u\rangle}}\right)^{-1}. 
\end{equation*} 
Then, the goal is to solve the maximize problem:
\begin{equation*}
 \max\limits_{\substack{\mathbf u\in\mathbb{R}^2\\\mathbf u\neq0}}{\frac{\vert\langle\mathbf u, A\mathbf u\rangle\vert}{\langle \mathbf u, B\mathbf u\rangle}}=\max\limits_{\substack{\mathbf u\in\mathbb{R}^2\\ \mathbf u\neq0}}{\frac{\vert\langle\mathbf u, B^{-\frac{1}{2}}AB^{-\frac{1}{2}}\mathbf u\rangle\vert}{\Vert \mathbf u\Vert^2}}.   
\end{equation*}
We know that $C=B^{-\frac{1}{2}}AB^{-\frac{1}{2}}$ is a symmetric matrix, and thus
\begin{equation*}
   \max\limits_{\substack{\mathbf u\in\mathbb{R}^2\\ \mathbf u\neq0}}{\frac{\vert\langle\mathbf u, B^{-\frac{1}{2}}AB^{-\frac{1}{2}}\mathbf u\rangle\vert}{\Vert \mathbf u\Vert^2}}=\rho(C).  
\end{equation*}
Therefore,
\begin{equation*}
    \min\limits_{\substack{\nu\in\mathbb{R}}}\frac{Q(\nu)}{\vert P(\nu)\vert}=\left(\max\limits_{\lambda \text{ eigenvalue of }B^{-1/2}AB^{-1/2}}\vert\lambda\vert\right)^{-1}.
\end{equation*}

\section*{Acknowledgments}
The work of Y. Liu was supported by UZH Postdoc Grant, 2024 / Verf\"{u}gung Nr. FK-24-110 and SNSF grant 200020$\_$204917.

\bibliographystyle{plain}
\bibliography{reference}

\begin{thebibliography}{10}

\bibitem{Abgrall_camc}
R.~Abgrall.
\newblock A combination of {R}esidual {D}istribution and the {A}ctive {F}lux
  formulations or a new class of schemes that can combine several writings of
  the same hyperbolic problem: application to the 1{D} {E}uler equation.
\newblock {\em Commun. Appl. Math. Comput.}, 5:370--402, 2023.

\bibitem{Abgrall_Lin_Liu}
R.~Abgrall, J.~Lin, and Y.~Liu.
\newblock Active {F}lux for triangular meshes for compressible flows problems.
\newblock {\em to appear in Beijing J. of Pure and Appl. Math.}, 2024.
\newblock arXiv:2312.11271.

\bibitem{AbgrallLiu}
R.~Abgrall and Y.~Liu.
\newblock A new approach for designing well-balanced schemes for the shallow
  water equations: a combination of conservative and primitive formulations.
\newblock {\em SIAM J. Sci. Comput.}, 46:A3375--A3400, 2024.

\bibitem{AF4}
W.~Barsukow.
\newblock The {A}ctive {F}lux scheme for nonlinear problems.
\newblock {\em J. Sci. Comput.}, 86, 2021.
\newblock Paper No. 3.

\bibitem{vonNeuman}
H.~A. Bethe, K.~Fuchs, J.~O. Hirschfelder, J.~L. Magree, R.~E. Peirls, and
  J.~von Neuman.
\newblock Blast waves.
\newblock Technical Report LA-2000, Los Alamos Scientific Laboratory, 1947.

\bibitem{Carlier_IDP}
V.~Carlier and F.~Renac.
\newblock Invariant domain preserving high-order spectral discontinuous
  approximations of hyperbolic systems.
\newblock {\em SIAM J. Sci. Comput.}, 45:A1385--A1412, 2023.

\bibitem{CHK_21}
E.~Chudzik, C.~Helzel, and D.~Kerkmann.
\newblock The {C}artesian grid {A}ctive {F}lux method: linear stability and
  bound preserving limiting.
\newblock {\em Appl. Math. Comput.}, 393:125501, 2021.

\bibitem{CDL}
S.~Clain, S.~Diot, and R.~Loub\`{e}re.
\newblock A high-order finite volume method for systems of conservation
  laws---multi-dimensional optimal order detection ({MOOD}).
\newblock {\em J. Comput. Phys.}, 230:4028--4050, 2011.

\bibitem{duan2024AF}
J.~Duan, W.~Barsukow, and C.~Klingenberg.
\newblock Active {F}lux methods for hyperbolic conservation laws -- flux vector
  splitting and bound-preservation.
\newblock {\em SIAM J. Sci. Comput.}, 47:A811--A837, 2025.

\bibitem{Dumbser_MOOD}
M.~Dumbser and R.~Loub\`{e}re.
\newblock A simple robust and accurate a posteriori sub-cell finite volume
  limiter for the discontinuous {G}alerkin method on unstructured meshes.
\newblock {\em J. Comput. Phys.}, 319:163--199, 2016.

\bibitem{Dumbser_MOOD2}
M.~Dumbser, O.~Zanotti, R.~Loub\`{e}re, and S.~Diot.
\newblock A posteriori subcell limiting of the discontinuous {G}alerkin finite
  element method for hyperbolic conservation laws.
\newblock {\em J. Comput. Phys.}, 278:47--75, 2014.

\bibitem{AF3}
T.~A. Eyman.
\newblock {\em Active {F}lux}.
\newblock PhD thesis, University of Michigan, 2013.

\bibitem{AF1}
T.~A. Eyman and P.~L. Roe.
\newblock Active {F}lux.
\newblock 49th AIAA Aerospace Science Meeting, 2011.

\bibitem{AF2}
T.~A. Eyman and P.~L. Roe.
\newblock Active {F}lux for systems.
\newblock 20 th AIAA Computationa Fluid Dynamics Conference, 2011.

\bibitem{Guermond_IDP}
J.-L. Guermond, M.~Nazarov, B.~Popov, and I.~Tomas.
\newblock Second-order invariant domain preserving approximation of the {E}uler
  equations using convex limiting.
\newblock {\em SIAM J. Sci. Comput.}, 40:A3211--A3239, 2018.

\bibitem{Guermond_IDP2}
J.-L. Guermond and B.~Popov.
\newblock Invariant domains and first-order continuous finite element
  approximation for hyperbolic systems.
\newblock {\em SIAM J. Numer. Anal.}, 54:2466--2489, 2016.

\bibitem{Guermond_IDP3}
J.-L. Guermond, B.~Popov, and I.~Tomas.
\newblock Invariant domain preserving discretization independent schemes and
  convex limiting for hyperbolic systems.
\newblock {\em Comput. Method. Appl. M.}, 347:143--175, 2019.

\bibitem{Hajduk_MCL}
H.~Hajduk.
\newblock Monolithic convex limiting in discontinuous {G}alerkin
  discretizations of hyperbolic conservation laws.
\newblock {\em Comput. Math. Appl.}, 87:120--138, 2021.

\bibitem{HKS}
C.~Helzel, D.~Kerkmann, and L.~Scandurra.
\newblock A new {ADER} method inspired by the {A}ctive {F}lux method.
\newblock {\em J. Sci. Comput.}, 80:35--61, 2019.

\bibitem{Jiang1996}
G.~S. Jiang and C.-W. Shu.
\newblock Efficient implementation of weighted {WENO} schemes.
\newblock {\em J. Comput. Phys.}, 126:202--228, 1996.

\bibitem{Kamm2000}
J.~R. Kamm.
\newblock Evaluation of the {S}edov-von {N}eumann-{T}aylor blast wave
  solutions.
\newblock Technical Report LA-UR-00-6055, Los Alamos National Laboratory, 2000.

\bibitem{Kuzmin_MCL1}
D.~Kuzmin.
\newblock Monolithic convex limiting for continuous finite element
  discretizations of hyperbolic conservation laws.
\newblock {\em Comput. Method. Appl. M.}, 361:112804, 2020.

\bibitem{Kuzmin_BP}
D.~Kuzmin, M.~Quezada~de Luna, D.~Ketcheson, and J.~Gr{\"u}ll.
\newblock Bound-preserving flux limiting for high-order explicit {R}unge
  {K}utta time discretizations of hyperbolic conservation laws.
\newblock {\em J. Sci. Comput.}, 91:21, 2022.

\bibitem{Liu_PAMPA_SW2D}
Y.~Liu.
\newblock Well-balanced {P}oint-{A}verage-{M}oment {P}olynomi{A}l-interpreted
  ({PAMPA}) methods for shallow water equations with horizontal temperature
  gradients on triangular meshes.
\newblock 2024.
\newblock arXiv:2409.12606.

\bibitem{LiuBarsukow}
Y.~Liu and W.~Barsukow.
\newblock An arbitrarily high-order fully well-balanced hybrid finite
  element-finite volume method for a one-dimensional blood flow model.
\newblock 2024.
\newblock arXiv:2404.18124.

\bibitem{Lohmann_FCT}
C.~Lohmann, D.~Kuzmin, J.~N. Shadid, and S.~Mabuza.
\newblock Flux-corrected transport algorithms for continuous {G}alerkin methods
  based on high order {B}ernstein finite elements.
\newblock {\em J. Comput. Phys.}, 344:151--186, 2017.

\bibitem{Loubere_MOOD}
R.~Loub\`{e}re, R.~Turpault, and A.~Bourriaud.
\newblock A {MOOD}-like compact high order finite volume scheme with adaptive
  mesh refinement.
\newblock {\em Appl. Math. Comput.}, 443:127792, 2023.

\bibitem{MaireHDR}
P.-H. Maire.
\newblock {\em Contribution to the numerical modeling of inertial confinement
  fusion}.
\newblock PhD thesis, University of Bordeaux, 2011.

\bibitem{Sedov}
L.~I. Sedov.
\newblock {\em {Similarity and dimensional methods in mechanics}}.
\newblock Academic Press, New York, 1959.

\bibitem{Shu1989}
C.-W. Shu and S.~Osher.
\newblock Efficient implementation of essentially non-oscillatory
  shock-capturing schemes, {II}.
\newblock {\em J. Comput. Phys.}, 83:32--78, 1989.

\bibitem{Vilar}
F.~Vilar.
\newblock A posterior correction of high-order discontinuous {G}alerkin scheme
  through subcell finite volume formulation and flux reconstruction.
\newblock {\em J. Comput. Phys.}, 387:245--279, 2019.

\bibitem{Vilar_DGFV}
F.~Vilar.
\newblock Local subcell monolithic {DG}/{FV} convex property preserving scheme
  on unstructured grids and entropy consideration.
\newblock {\em J. Comput. Phys.}, 521:113535, 2025.

\bibitem{wu2023geometric}
K.~Wu and C.-W. Shu.
\newblock Geometric quasilinearization framework for analysis and design of
  bound-preserving schemes.
\newblock {\em SIAM Rev.}, 65:1031--1073, 2023.

\bibitem{Zhang2010}
X.~Zhang and C.-W. Shu.
\newblock On maximum-principle-satisfying high order schemes for scalar
  conservation laws.
\newblock {\em J. Comput. Phys.}, 229:3091--3120, 2010.

\bibitem{Zhang_MP}
X.~Zhang, Y.~Xia, and C.-W. Shu.
\newblock Maximum-principle-satisfying and positivity-preserving high order
  discontinuous {G}alerkin schemes for conservation laws on triangular meshes.
\newblock {\em J. Sci. Comput.}, 50:29--62, 2012.

\end{thebibliography}

\end{document}